\documentclass[12pt]{elsarticle}

\usepackage{amsmath}
\usepackage{amssymb}
\usepackage{graphicx}
\usepackage{hyperref}
\usepackage{xcolor}
\usepackage{multicol}
\usepackage{calrsfs}
\usepackage{algorithm}
\usepackage{algorithmic}
\usepackage{cancel}
\usepackage{multirow}
\usepackage[algo2e]{algorithm2e} 
\newcommand{\RR}{{\mathbb R}}

\newcommand{\NN}{\mathbb{N}}

\newcommand{\C}{\mathcal{C}}
\newcommand{\V}{\mathcal{V}}
\newcommand{\W}{\mathcal{W}}
\newcommand{\bigO}{\mathcal{O}}

\newcommand\MVB[1]{{\color{black}#1}}
\newcommand\WT[1]{{\color{black}#1}}
\newcommand\MC[1]{{\color{black}#1}}

\begin{document}
	
	\title{A Parametric Family of  Polynomial Wavelets for  Signal and Image Processing}
	
	\author[inst1]{Mariantonia Cotronei\corref{cor1}}
    \ead{mariantonia.cotronei@unirc.it}
	\author[inst2]{Woula Themistoclakis}
    \ead{woula.themistoclakis@cnr.it}
	\author[inst3]{Marc Van Barel}
    \ead{marc.vanbarel@kuleuven.be}
	
	\affiliation[inst1]{organization={DIIES, Università Mediterranea di Reggio Calabria},
		addressline={Via Zehender Loc. Feo di Vito},
		city={Reggio Calabria},
		postcode={89122},
		country={Italy}}
\cortext[cor1]{Corresponding author}

    \affiliation[inst2]{organization={C.N.R. National
        Research Council of Italy,
        IAC Institute for Applied Computation ``Mauro Picone''},
		addressline={Via P. Castellino 111},
		city={Naples},
		postcode={80131},
		country={Italy}}
		
	\affiliation[inst3]{organization={Department of Computer Science, KU Leuven},
			addressline={Celestijnenlaan 200A},
			city={Heverlee},
			postcode={B-3001},
			country={Belgium}}
	
	\begin{abstract}
	This paper investigates the potential applications of a parametric family of polynomial wavelets that has been recently introduced starting from de la Vall\'ee Poussin (VP) interpolation at Chebyshev nodes. Unlike classical wavelets, which are constructed on the real line, these VP wavelets are defined on a bounded interval, offering the advantage of handling boundaries naturally while maintaining computational efficiency. In addition, the structure of these wavelets enables the use of fast algorithms for decomposition and reconstruction. Furthermore,  the flexibility offered by a free parameter allows a better control of localized singularities, such as edges in images. On the basis of previous theoretical foundations, we show the effectiveness of the VP wavelets for basic signal denoising and image compression, emphasizing their potential for more advanced signal and image processing tasks.
	\end{abstract}
	
	\begin{keyword}
		Wavelets; Multiresolution analysis; De la Vall\'ee Poussin  interpolation; Sparse representations; Denoising; Compression
	\end{keyword}
	
	\maketitle
	
	\section{Introduction}
Wavelet transforms are widely recognized as powerful tools in various applications, particularly in signal and image processing, being able to effectively extract frequency content from data and represent them using a reduced set of information.

Our goal is to contribute to the extensive research developed over the past decades by implementing and testing a special class of wavelet transforms based on a non-standard multiresolution analysis in which the resolution level is tripled instead of doubled at each step. Specifically, we consider a parametric family of \WT{polynomial} wavelets recently introduced in \cite{TheVanBar2024}. \WT{Similarly to other kinds of polynomial wavelets (see, e.g. \cite{FischerPrestin1997, FT2002,KilgorePrestin1996}) there is 
no corresponding scaling or wavelet mother functions neither any dilation and translation is implemented, but, at any resolution level, we have different, well localized, polynomial scaling and wavelet basis functions generating the appproximation and detail spaces, respectively. }

\MC{The polynomial wavelets examined in this paper are based on de la Vallée Poussin (VP) interpolation at Chebyshev nodes \cite{Th-1999, Th-2012, Th-L1}. These VP wavelets are interpolating, non-orthogonal polynomials defined on the interval $[-1,1]$, to which any other compact interval can be mapped. They can also be naturally extended to higher dimensions through tensor-product constructions. Hence, data defined on finite domains can be processed directly by VP wavelets, without requiring additional modifications, unlike classical wavelets, which are intrinsically defined on the entire real line through translations and dilations of a single mother function.}.

 Furthermore, even with the uniform norm, they  possess good approximation properties, which can be further controlled through a free parameter.
 
In addition, the structure of the VP wavelets allows for the implementation of fast algorithms for data decomposition and reconstruction, ensuring computational efficiency. 

In this paper, after presenting the main properties of the new wavelet basis, we describe the implementation of the corresponding transform and inverse transform, providing details relevant to their correct application in signal processing.

We then focus on two classical applications: signal denoising and image compression. In such contexts, our goal is not to develop new algorithms, but to demonstrate the potential performance of the proposed basis, particularly in comparison with standard wavelet filters.

The paper is organized as follows: Section 2 introduces the univariate VP scaling and wavelet functions, detailing their construction and key properties. Section 3 presents the decomposition and reconstruction algorithms, along with details of a correct application in signal processing. Sections 5 and 6 discuss the application of VP wavelets to signal denoising and image compression, respectively, providing experimental results and comparisons with traditional wavelet methods. Finally, conclusions are drawn, highlighting potential future directions.

\section{Univariate VP scaling and wavelet functions}
For any $n\in\NN$, we denote the zeros of the Chebyshev polynomial $T_n(x)=\cos(n\arccos{x})$ by
\[
x_k^n:=\cos\frac{(2k-1)\pi}{2n},\qquad k=1,\ldots,n.
\]
They are also zeros of $T_{3n}(x)$ and specifically we have
\[
x_k^n=x_{3k-1}^{3n}, \qquad k=1,\ldots, n.
\]
The remaining zeros of $T_{3n}$, which are not zeros of $T_n$, are denoted by $y_k^n$ with $k=1,\ldots,2n$. More precisely, we set
\[
y_{2k-1}^n:=x_{3k-2}^{3n} \quad\mbox{and}\quad y_{2k}^n:=x_{3k}^{3n},\qquad k=1,\ldots,n.
\]
The two sets of Chebyshev zeros
\[ 
X_n:=\{x_k^n : \ k=1,\ldots,n\}\quad\mbox{and}\quad
Y_n=\{y_k^n : \ k=1,\ldots, 2n\}
\]
constitute the interpolation nodes of the scaling and wavelet functions, which we are going to define below.

To this aim, we make use of the following polynomial kernel \cite{TheVanBar2024}
\begin{equation}\label{VP-kernel}
\phi_{n}^m(x,y)=\frac 2n\sum_{r=0}^{n+m-1}\hspace{-.1cm}{}^\prime
\mu_{n,r}^m T_r(x)T_r(y),\qquad n,m\in\NN, \quad 0<m<n,  \end{equation}
where the notation $\sum {}^\prime$ means that the first addendum of the summation is halved, and the coefficients $\mu_{n,r}^m$ are defined as follows
\vspace{.1cm}\newline
\begin{multicols}{2}
 $\mu_{n,r}^m :=\left\{\begin{array}{ll}
1 & \mbox{if}\ 0\le r\le n-m, \\ [.1in]
\displaystyle\frac{m+n-r}{2m} & \mbox{if}\
n-m< r< n+m, 
\end{array}\right.$ 
\columnbreak
 
\hspace{1.5cm}
\setlength{\unitlength}{.85cm}
\begin{picture}(7,2)
\put(0.3,0.3){\circle*{.125}}\put(0,0){\makebox(0,0){\scriptsize 0}}
\put(0.3,0.3){\vector(1,0){5}}\put(5.5,0.3){\makebox(0,0){\footnotesize r}}
\put(0.3,0.3){\vector(0,1){2.5}}
\put(3.5,1.5){\makebox(0,0){\color{blue}\footnotesize$\mu_{n,r}^m$}}
\put(0,2){\makebox(0,0){ \scriptsize 1}}\put(0.3,2){\circle*{.125}}
\put(0.3,2){\color{blue}\line(1,0){1.7}}
\put(2,2){\color{blue}\line(1,-1){1.7}}
\put(2,0){\makebox(0,0){\footnotesize n-m}}\put(2,0.3){\circle*{.125}}
\put(2,0.6){\circle*{.08}}\put(2,0.9){\circle*{.08}}
\put(2,1.2){\circle*{.08}}\put(2,1.5){\circle*{.08}}\put(2,1.8){\circle*{.08}}
\put(3.7,0){\makebox(0,0){  \footnotesize n+m}}\put(3.7,0.3){\circle*{.125}}
\end{picture}
\end{multicols}
Equivalently, $\phi_n^m$ can be written in the following trigonometric form \cite{OT-APNUM21}
\[
\phi_n^m(\cos t, \cos s)=\frac 1{4 nm}\left(\frac{\sin[n(t-s)]\sin[m(t-s)]}{\sin^2[(t-s)/2]}
+ \frac{\sin[n(t+s)]\sin[m(t+s)]}{\sin^2[(t+s)/2]}\right) .
\]

Under the previous notation, for all $x\in[-1,1]$, each resolution level $n\in\NN$ and any integer parameter $0<m<n$, the {\em VP scaling functions} 
are defined as follows \cite{TheVanBar2024}
\begin{equation}\label{sca}
\phi_{n,k}^m(x)=\phi_n^m(x_k^n,\ x),\qquad k=1,\ldots,n,
\end{equation}
while the {\em VP wavelets} 
are 
\begin{equation}\label{wav}
\psi_{n,k}^m(x)=\phi_{3n}^m(y_k^n,x) - \sum_{h=1}^n \phi_{n,h}^m(y_k^{n})\phi_{3n}^m(x_h^n,x), \qquad k=1,\ldots,2n.
\end{equation}

\MC{Note that the terms “scaling function” and “wavelet” are used here in a broader sense. Although these functions are not generated by dyadic dilations and translations of a single prototype function, they share several key wavelet-like features, as illustrated below. This is in line with earlier polynomial constructions already called "wavelets" (e.g.,~\cite{FischerPrestin1997,FT2002,KilgorePrestin1996}).}

Let us summarize the main properties of these polynomials valid for any pair of integers $0<m<n$ \cite{TheVanBar2024}.
\begin{enumerate}
    \item {\it Interpolation}. Using the Kronecker delta, we have
\begin{eqnarray}
  \phi_{n,k}^m(x_h^n)&=& \delta_{h,k},   \qquad h,k=1,\ldots, n,  \\
   \psi_{n,k}^m(y_h^n)&=& \delta_{h,k},   \qquad h,k=1,\ldots, 2n. 
\end{eqnarray}    
\item {\it Moments}. Setting $dw(x)=\frac{dx}{\sqrt{1-x^2}}$, for all $r=0,1,\ldots,n-m$, we have
\begin{eqnarray}
\int_{-1}^1 x^r \phi_{n,k}^m(x)dw(x)&=& (x_k^n)^r, \qquad k=1,\ldots,n,\\ 
\int_{-1}^1 x^r \psi_{n,k}^m(x)dw(x)&=& 0, \qquad\qquad k=1,\ldots, 2n.
\end{eqnarray}
\item {\it Stability}. Setting, for any function $f$,
\[ 
\|f\|_{p} =\left\{
\begin{array}{ll}
\displaystyle \left(\int_{-1}^1|f(x)|^p w(x)dx\right)^\frac 1p & 
1\le p<\infty,\\ [.1in]
\displaystyle \sup_{|x|\le 1}|f(x)| & p=\infty,\end{array}\right.
\]
and, for any vector $\vec{v}_N=(v_1,\ldots,v_N)\in\RR^N$,
\[
\|\vec{v}_N\|_{\ell^p}=\left\{\begin{array}{ll}
\displaystyle \left(\frac \pi{N}\sum_{k=1}^N |v_{k}|^p\right)^\frac 1p & 1\le p<\infty,\\
\displaystyle \max_{1\le k\le N}|v_{k}| & p=\infty,
\end{array}\right. 
\]
we have that for any $1\le p\le\infty$ and all vectors $\vec{a}_{n}=(a_{1}, \ldots, a_{n})\in\RR^{n}$, and $\vec{b}_{2n}=(b_{1}, \ldots, b_{2n})\in\RR^{2n}$, the following inequalities
\begin{eqnarray}
\label{Riesz-sca}
\C_1 \|\vec{a}_{n}\|_{\ell^p}\le  & \left\|\sum_{k=1}^{n} a_{k}\phi_{n,k}^m\right\|_p\le &\C_2\left(1+\log\frac nm\right)\ \|\vec{a}_{n}\|_{\ell^p}, \\
\label{Riesz-wav}
\C_3 \|\vec{b}_{2n}\|_{\ell^p}\le  & \left\|\sum_{k=1}^{2n} b_{k}\psi_{n,k}^m\right\|_p\le &\C_4\left(1+\log^2\frac nm\right)\ \|\vec{b}_{2n}\|_{\ell^p}, 
\end{eqnarray}
hold, where $\C_i>0$, $i=1,..,4$, are constants independent of $\vec{a}_n$, $\vec{b}_{2n}$, $n$ and $m$.
\item {\it Convergence}.  For any fixed $\theta\in ]0,1[$, if we take $m=\lfloor \theta n\rfloor$ and a sufficiently large resolution level $n\in\NN$, then any continuous function $f$ on $[-1,1]$ can be uniformly approximated with the desired precision by the VP interpolation polynomial \cite{Th-2012}
\begin{equation}\label{Vnm}
V_n^m f(x):=\sum_{k=1}^nf(x_k^n)\phi_{n,k}^m(x),\qquad  x\in [-1,1],
\end{equation}
namely, we have
\begin{equation}\label{lim}
\lim_{\substack{
n\to\infty\\
m=\lfloor \theta n\rfloor}}
\|f-V_n^mf\|_\infty=0, \qquad \forall f\in C^0[-1,1],\qquad \forall\theta\in ]0,1[ .
\end{equation}
Moreover, we have the same convergence order of the error of best polynomial approximation $E_n(f):=\inf_{\deg (P)\le n}\|f-P\|_\infty$, which depends on the smoothness degree of $f$ \cite{DT}. In particular, if $f$ is H\"older continuous with exponent $0<\alpha\le 1$ then we have
\begin{equation}\label{V-lip}
\|f-V_n^mf\|_\infty=\bigO(n^{-\alpha})
\end{equation}
Furthermore, if $f$ is continuously differentiable up to order $2s$, for some $s\in\NN$, then we have \cite[Thm. 5.3]{OT-DRNA21} 
\begin{equation}\label{V-der}
 \displaystyle \|f^{(r)}-(V_n^mf)^{(r)}\|_\infty=\bigO(n^{-2s+r}), \qquad r=0,1,\ldots,s
\end{equation}
where $f^{(r)}$ denotes the $r$-th derivative of $f$, and $f^{(0)}:=f$.
\end{enumerate}

\section{Decomposition and reconstruction algorithms}\label{sec:algo}
\MC{The above defined VP scaling and wavelet functions generate a hierarchical structure of approximation and detail spaces. 
At a given resolution level $n \in \NN$ and for an arbitrarily fixed parameter $m<n$, 
we define
\begin{eqnarray}
\label{Vn}
 \V_n^m&=&\operatorname{span}\, \{\phi_{n,k}^m, \ k=1,\ldots,n\}, \qquad  \dim \V_n^m=n, \\
 \label{Wn}
 \W_n^m&=& \operatorname{span}\, \{\psi_{n,k}^m, \ k=1,\ldots,2n\}, \qquad  \dim \W_n^m=2n.
\end{eqnarray}
Each space $\V_n^m$ is finite-dimensional and therefore a closed subspace of the weighted
Hilbert space $L^2_w([-1,1])$, endowed with the scalar product
$$
\langle f,g\rangle := \int_{-1}^1 f(x)g(x)\,dw(x), 
\qquad dw(x)=\frac{dx}{\sqrt{1-x^2}} .
$$
Moreover, these spaces satisfy \cite{TheVanBar2024}
\begin{equation}\label{v+w}
\V^m_{3n}=\V^m_{n}\oplus \W_n^m, 
\qquad 
\V_n^m\perp\W_n^m, 
\qquad \forall\, n>m\in\NN.
\end{equation}

\noindent\textbf{Remark.}
The relation \eqref{v+w} defines a nested sequence of polynomial approximation and detail
spaces that mirrors the hierarchical refinement of a multiresolution
analysis (MRA). However, since no translation or dilation operators are involved and $\V_n^m$ is not a shift-invariant subspace of $L^2(\mathbb{R})$, this construction does not yield a
"classical" \emph{multiresolution analysis}  (MRA). In this paper, this term is therefore used in a
nonstandard or generalized sense to highlight the hierarchical and orthogonal
decomposition \eqref{v+w}
which refines the resolution by a factor of three on the compact interval $[-1,1]$.
}

\bigskip

\MC{In virtue of \eqref{v+w},} any function $f_{3n}\in \V_{3n}^m$ can be uniquely decomposed as
\begin{equation}\label{f+g}
f_{3n}=f_n+g_{2n}, \qquad \mbox{with $\quad f_n\in\V_n^m\quad$ and 
$\quad g_{2n}\in\W_n^m$}.
\end{equation} 
Using the interpolating bases \eqref{Vn}--\eqref{Wn}, these functions are uniquely determined by their basis coefficients, i.e., we have
\begin{eqnarray}\label{fn}
f_n(x)&=&\sum_{k=1}^n a_{n,k}\phi_{n,k}^m(x),\qquad a_{n,k}=f_n(x_k^n),\qquad \forall f_n\in\V_n^m\\
\label{gn}
g_{2n}(x)&=&\sum_{k=1}^{2n} b_{n,k}\psi_{n,k}^m(x),\qquad b_{n,k}=g_{2n}(y_k^n),\qquad \forall g_{2n}\in\W_n^m .
\end{eqnarray}
Fast decomposition and reconstruction algorithms can be  given in order to compute the basis coefficients of $f_n$ and $g_{2n}$ from the basis coefficients of $f_{3n}$, and viceversa \cite[Thm. 5.2]{TheVanBar2024}. 

Setting
\[
p_r(x)=\sqrt{\frac 2\pi}\frac{T_{r-1}(x)}{\sqrt{1+\delta_{r,1}}}, \qquad r=1,2,\ldots
\]
and 
\[
q_r(x)=\left\{\begin{array}{lr}
  p_r(x)   & r=1,\ldots, n-m+1  \\
  \frac{m+(n-r+1)}{2m}p_r(x)- \frac{m-(n-r+1)}{2m}p_{2n-r+2}(x)  & 
  r=n-m+2,\ldots, n
\end{array}\right.
\]
these algorithms are based on the following discrete transforms
\begin{itemize}
    \item $\vec{u}=\mathcal{T}_1(\vec{v})$ that transforms $\vec{v}=(v_1,\ldots, v_n)$ in $\vec{u}\in\RR^n$  with entries
    \[
    u_r=\sum_{s=1}^n v_s p_r(x_s^n),\qquad r=1,\ldots, n
    \]
\item $\vec{u}=\mathcal{T}_2(\vec{v})$ that transforms $\vec{v}=(v_1,\ldots, v_n)$ in $\vec{u}\in\RR^n$  with entries
    \[
    u_s=\sum_{r=1}^n v_r p_r(x_s^n),\qquad s=1,\ldots, n
    \]
    \item $\vec{u}=\mathcal{T}_3(\vec{v})$ that transforms $\vec{v}=(v_1,\ldots, v_{2n})$ in $\vec{u}\in\RR^n$  with entries
    \[
    u_r=\sum_{s=1}^{2n} v_s q_r(y_s^n),\qquad r=1,\ldots, n
    \]
    \item $\vec{u}=\mathcal{T}_4(\vec{v})$ that transforms $\vec{v}=(v_1,\ldots, v_n)$ in $\vec{u}\in\RR^{2n}$  with entries
    \[
    u_s=\sum_{r=1}^n v_r q_r(y_s^n),\qquad s=1,\ldots, 2n
    \]
\end{itemize}
By using the previous transforms the decomposition and reconstruction algorithms can be schematized as follows

\begin{minipage}[t]{0.5\textwidth}
\begin{algorithm}[H]
\caption{Decomposition \WT{(classic)}}
\KwIn{$\vec{a}_{3n}$
 basis coeff. of $f_{3n}=f_n+g_{2n}$}
\KwOut{$\vec{a}_n$ and $\vec{b}_n$ basis coeff.of $f_n$ and $g_{2n}$}\vspace{.2cm}
\begin{algorithmic}[1] 
\FOR{$k=1:n$}\vspace{.1cm}
\STATE $v_k=a_{3n,3k-1}$
\STATE $u_{2k-1}=a_{3n,3k-2}$
\STATE $u_{2k}=a_{3n,3k}$\vspace{.1cm}
\ENDFOR \vspace{.2cm}
\STATE $\displaystyle \vec{x}=\frac\pi{3n}\mathcal{T}_1\left( \vec{v}\right)$\vspace{.1cm}
\STATE $\displaystyle\vec{y}=\frac\pi{3n}\mathcal{T}_3\left( \vec{u}\right)$
\vspace{.2cm}
\FOR{$k=1:n-m+1$}\vspace{.1cm}
\STATE $v_k=x_k+y_k$\vspace{.1cm}
\ENDFOR\vspace{.1cm}
\FOR{$k=n-m+2:n$}\vspace{.2cm}
\STATE $v_k=\frac{2m^2(x_k+y_k)}{m^2+(n-k+1)^2}$\vspace{.2cm}
\ENDFOR
\STATE $\vec{a}_n=\mathcal{T}_2(\vec{v})$
\STATE $\vec{b}_n=\mathcal{T}_4(\vec{v})$
\end{algorithmic}
\end{algorithm}
\end{minipage}%
\hspace{.15cm}
\begin{minipage}[t]{0.5\textwidth}
\begin{algorithm}[H]
\caption{Reconstruction \WT{(classic)}}
\KwIn{$\vec{a}_n$ and $\vec{b}_n$ basis coeff.of $f_n$ and $g_{2n}$}
\KwOut{$\vec{a}_{3n}$
 basis coeff. of $f_{3n}=f_n+g_{2n}$}\vspace{.2cm}
\begin{algorithmic}[1]
 \STATE $\displaystyle \vec{x}=\frac\pi{n}\mathcal{T}_1\left( \vec{a}_n\right)$\vspace{.2cm}
\STATE $\displaystyle\vec{y}=\frac\pi{3n}\mathcal{T}_3\left( \vec{b}_n\right)$\vspace{.2cm}
\STATE $\vec{u}= \vec{a}_n-\mathcal{T}_2(\vec{y})$\vspace{.2cm}
\STATE $\vec{v}= \vec{b}_n+ \mathcal{T}_4(\vec{x})$\vspace{.2cm}
\FOR{$k=1:n$}\vspace{.1cm}
\STATE $a_{3n,3k-1}=u_k$
\STATE $a_{3n,3k-2}=v_{2k-1}$
\STATE $a_{3n,3k}=v_{2k}$\vspace{.1cm}
\ENDFOR
\end{algorithmic}
\end{algorithm}
\end{minipage}
\smallskip

As observed in \cite{TheVanBar2024}, each transform $\mathcal{T}_i$ can be efficiently computed using 
the well--known {\em discrete cosine transform} (DCT) and its inverse (IDCT), in Matlab known as DCT of type 2 and 3, respectively. They transform any $\vec{v}=(v_1,\ldots,v_n)\in\RR^n$ in $\vec{u}\in\RR^n$ defined as follows
\begin{itemize}
\item
$\vec{u}=DCT (\vec{v})$ has entries
$\displaystyle
u_r=\sqrt{\frac 2n}\frac 1{\sqrt{1+\delta_{r,1}}}\sum_{s=1}^n
v_s\ T_{r-1}(x_s^n),
\quad r=1,\ldots,n$
\item
$\vec{u}=IDCT (\vec{v})$ has entries
$\displaystyle 
u_s=\sqrt{\frac 2n}\sum_{r=1}^n \frac {v_r}{\sqrt{1+\delta_{r,1}}}
\ T_{r-1}(x_s^n),
\quad s =1,\ldots,n$.
\end{itemize}
Using these transforms, we have
\[
\mathcal{T}_1(\vec{v})=\sqrt{\frac n\pi} DCT(\vec{v}) \quad\mbox{and}\quad\mathcal{T}_2(\vec{v})=\sqrt{\frac n\pi} IDCT(\vec{v}), \qquad \forall \vec{v}\in\RR^n.
\]
Moreover, the transforms $\mathcal{T}_3$ and $\mathcal{T}_4$ can be computed using the following algorithms.
\begin{minipage}[t]{0.55\textwidth}
\begin{algorithm}[H]
\caption{Compute $\vec{u}=\mathcal{T}_3(\vec{v})$}
{\small
\KwIn{$\vec{v}=(v_1,\ldots, v_{2n})$}
\KwOut{$\vec{u}=(u_1,\ldots,u_{n})$}\vspace{.2cm}
}
\begin{algorithmic}[1] 
\FOR{$k=1:n$}\vspace{.1cm}
\STATE $w_{3k-2}=v_{2k-1}$
\STATE $w_{3k-1}=0$
\STATE $w_{3k}=v_{2k}$\vspace{.1cm}
\ENDFOR \vspace{.2cm}
\STATE $
\vec{x}=\sqrt{\frac{3n}\pi} DCT(\vec{w})$
\vspace{.2cm}
\FOR{$k=1:n-m+1$}\vspace{.1cm}
\STATE $u_k=x_k$\vspace{.1cm}
\ENDFOR\vspace{.1cm}
\FOR{$k=n-m+2:n$}\vspace{.2cm}
\STATE $
\hspace{-.2cm} u_k=\frac{m+(n-k+1)}{2m} x_k - \frac{m-(n-k+1)}{2m}x_{2n-k+2}$\vspace{.2cm}
\ENDFOR
\end{algorithmic}
\end{algorithm}
\end{minipage}%
\hspace{.2cm}
\begin{minipage}[t]{0.45\textwidth}
\begin{algorithm}[H]
{\small
\caption{Compute $\vec{u}=\mathcal{T}_4(\vec{v})$}
\KwIn{$\vec{v}=(v_1,\ldots, v_n)$}
\KwOut{$\vec{u}=(u_1,\ldots,u_{2n})$}\vspace{.2cm}
\begin{algorithmic}[1]
{\footnotesize
\FOR{$k=1:n-m+1$}
\STATE $w_k=v_k$
\ENDFOR
\FOR{$k=n-m+2:n$} \vspace{.2cm}
\STATE $
w_k=\frac{m+(n-k+1)}{2m} v_k$\vspace{.1cm}
\STATE $
w_{2n-k+2}=- \frac{m-(n-k+1)}{2m}v_k$\vspace{.1cm}
\ENDFOR
\FOR{$k=n+m+4:3n$}
\STATE $w_k=0$
\ENDFOR \vspace{.2cm}
\STATE $\displaystyle
\vec{x}=\sqrt{\frac{3n}\pi} IDCT(\vec{w})$\vspace{.1cm}
\FOR{$k=1:n$}
\STATE $u_{2k-1}=x_{3k-2}$
\STATE $u_{2k}=x_{3k}$
\ENDFOR
}
\end{algorithmic}
}
\end{algorithm}
\end{minipage}
\bigskip

\MC{We conclude the section with some remarks about the computational complexity. 
The proposed  transform requires $ O(N \log N)$ operations per level due to its DCT-based implementation, compared with $ O(N)$ for standard dyadic wavelets. However, its  ternary decomposition structure
($\V_{3n}^{m} = \V_{n}^{m} \oplus \W_{n}^{m}$) allows more substantial data reduction per step, resulting in fewer decomposition levels and an overall computational cost that remains practically competitive while maintaining reconstruction quality, as confirmed in the following experimental sections.}

\WT{
\section{Application to signal processing}
}
The application of VP bases to signal processing requires some considerations.

First of all, let us discuss the mathematical model to adopt. As it is well known, a digital signal is mathematically represented by a vector, say $\vec{v}=(v_1,\ldots, v_N)$, of a certain length $N$. Its elements are typically considered as the values assumed by a certain function on $N$ equidistant nodes of an interval, which we can assume to be $[-1,1]$ without violating generality. In our model, we consider equidistant nodes not on the segment $[-1,1]$ but on the semicircle, i.e. 
\[
t_k=\frac{(2k-1)\pi}{2N}, \qquad k=1,\ldots,N,
\]
and we assume that behind the signal there exists a function $f$, representing it in the continuum, such that
\begin{equation}\label{eq.sample}
v_k= f(\cos t_k), \qquad k=1,\ldots,N.   
\end{equation}
Thus, using these samples we can approximate $f$, and hence the signal, by its  VP interpolation polynomial 
\begin{equation}\label{eq.VP}
V_N^mf(x)=\sum_{k=1}^Nv_k \phi_{N,k}^m(x), \quad m=\lfloor\theta N\rfloor, \quad \theta\in ]0,1[, \qquad x\in [-1,1],
\end{equation}
\WT{which represents the interpolating projection onto the approximation space $\V_N^m$, for arbitrary values of the parameter $\theta$, and hence $m$.

The non-uniform sampling \eqref{eq.sample} at the Chebyshev zeros and the consequent VP approximation \eqref{eq.VP} of a digital signal $\vec{v}=(v_1,\ldots, v_N)$ represent the non standard mathematical model adopted here. We recall that this model has been recently applied in \cite{Image_VP2D,Image_VP3D} where, taking the multivariate extension via tensor product, the authors get efficient and competitive algorithms for zooming arbitrary 2D and 3D images just resampling the multivariate VP approximation at the Chebyshev grid of the desired size. 

In this paper, we start from the same sampling/approximation model with the aim of handling
the signal $\vec{v}$,   interpreted as a scaling vector at the resolution $N$, with the usual wavelet techniques: the decomposition and reconstruction algorithms, respectively described in Algorithm~1 and Algorithm~2. 

To this aim, now, we have to face} two problems: the first concerns $N$, which according to our scheme should be a multiple of three, the second concerns the conservation of the signal energy.

\MC{To address the first issue, before applying the decomposition step described in Algorithm 1, we enlarge the input scaling vector $\vec{v}$ whenever its length is not a multiple of three, by adding up to two extra elements so that its size becomes divisible by three. 

Different strategies can be used to define these added elements. In Algorithm 5, we introduce a parameter {\it ind} controlling the extension rule: if {\it ind}$=0$, the added elements are set to 0; if \textit{ind}$=1$, they replicate the last entry $v_N$ of $\vec{v}$. Since the VP approximation preserves constants, to avoid introducing spurious oscillations we fixed {\it ind}$=1$ in all the experiments reported in this paper. However, when we set {\it ind}$=0$, the numerical results were essentially unchanged. In fact, because at most two elements are added before each decomposition step, the effect of the chosen rule becomes negligible, especially for large input sizes. 

During reconstruction, the previously added elements are removed to restore the original signal length. This truncation is performed level by level: before each reconstruction step, we check whether the length of the scaling vector is twice that of the corresponding detail vector. If not, we deduce that the scaling vector had been enlarged during decomposition, and we discard its last added entries so that the length ratio is restored. After the final reconstruction, the recovered signal may differ from the original by at most two samples, which corresponds to a very small information loss. In practice, the effect of this truncation on both the numerical and visual results was found to be insignificant.}

\MC{Regarding the second problem, since the VP scaling and wavelet functions are interpolating and not orthogonal, the Parseval identity is not applicable, and appropriate correction factors must be considered both in decomposition and reconstruction in order to assure energy conservation. In the classical wavelet setting, the scaling and wavelet coefficients are usually multiplied by $\sqrt{2}$ in the decomposition step, and divided by $\sqrt{2}$ in the reconstruction step.} In our case, we  multiply/divide the scaling coefficients by $\sqrt{3}$ and the wavelet coefficients by $\sqrt{3/2}$. To justify this, we observe that $f_{3n}=f_n+g_{2n}$ and $f_n\bot g_{2n}$ imply
{\begin{equation}\label{f+g-norm}
\|f_{3n}\|_2^2=\|f_n\|_2^2+\|g_{2n}\|_2^2.
\end{equation}}

On the other hand, by \eqref{fn}-\eqref{gn}, and 
\eqref{Riesz-sca}--\eqref{Riesz-wav} we get
\[
\|f_{3n}\|_2^2\sim \left(\frac \pi{3n}\sum_{k=1}^{3n}a_{3n,k}^2\right),\quad 
\|f_{n}\|_2^2\sim \left(\frac \pi{n}\sum_{k=1}^{n}a_{n,k}^2\right),\quad
\|g_{2n}\|_2^2\sim \left(\frac \pi{2n}\sum_{k=1}^{2n}b_{n,k}^2\right),
\]

where by $A\sim B$ we mean that there exist absolute constants $\C_1,\C_2>0$ (independent of the resolution level and the coefficients) such that $\C_1 A\le B\le \C_2 A$.

Consequently, from \eqref{f+g-norm} we deduce
{\[
 {\frac \pi{3n}} \sum_{k=1}^{3n}a_{3n,k}^2 
\sim 
{\frac\pi{n}}\sum_{k=1}^{n}a_{n,k}^2 +
 {\frac \pi{2n}} \sum_{k=1}^{2n}b_{n,k}^2 
\]}
that yields
{\begin{equation}\label{energy}
 \sum_{k=1}^{3n}\,a_{3n,k}^2 \sim
 \sum_{k=1}^{n}(\sqrt{3}a_{n,k})^2+
 \sum_{k=1}^{2n}\left(\sqrt{\frac 3 2}\,b_{n,k}\right)^2  .
\end{equation}}

\WT{
With all the above observations in mind, in the following Algorithm 5 and Algorithm 6, we summarize all the adjustments to the decomposition and reconstruction Algorithms 1 and 2, respectively. As we have previously discussed, these algorithms have no assumption on the size, which can also be not a multiple of three. Moreover, we fix by default $ind=1$ regarding the method to enlarge the scaling vectors, and we use the normalization factors $f_{sca}= \sqrt{3}$ and  $f_{wav}=\sqrt{3/2}$ to multiply (in decomposition) and divide (in reconstruction) with the scaling and detail vectors, respectively.
\begin{algorithm}[H]
\caption{Decomposition (any size)}
\KwIn{$\vec{v}=(v_1,\ldots, v_{N})$
scaling vector to decompose}
\KwOut{$\vec{u}$ scaling vector \\
\hspace{1.5cm} $\ \vec{w}$ detail vector}\vspace{.2cm}
\begin{algorithmic}[1] 
\STATE $ind=1$; $f_{sca}=\sqrt{3}$; $f_{wav}=\sqrt{\frac 32}$
\STATE $R =$ remainder term  of $N/3$
\IF{$R \ne 0$}
\FOR{$k=1: 3-R$}
\IF{$ind=1$}
    \STATE $v_{N+k}=v_N$
\ELSE
    \STATE $v_{N+k}=0$
\ENDIF
\ENDFOR
\ENDIF
\STATE Run Algorithm 1 with input vector $\vec{a}_{3n}=\vec{v}$ obtaining as output the scaling vector $\vec{a}_n$ and the detail vector $\vec{b}_{2n}$ 
\STATE $\vec{u}= f_{sca}\cdot \vec{a}_n$
\STATE $\vec{w}= f_{wav} \cdot \vec{b}_{2n}$
\end{algorithmic}
\end{algorithm}

\begin{algorithm}[H]
\caption{Reconstruction (any size)}
\KwIn{$\vec{u}$ scaling vector \\
\hspace{1.5cm} $\ \vec{w}$ detail vector}
\KwOut{$\vec{v}=(v_1,\ldots, v_{N})$
reconstructed scaling vector}\vspace{.2cm}
\begin{algorithmic}[1] 
\STATE $f_{sca}=\sqrt{3}$; $f_{wav}=\sqrt{\frac 32}$
\STATE $n=length(\vec{w})/2$
\STATE $\vec{a}_n = \vec{u}(1:n)$; 
\STATE $\vec{a}_n = \vec{a}_n/fsca$; 
\STATE $\vec{b}_{2n}=\vec{w}/ fwav$
\STATE Run Algorithm 2 with input vectors $\vec{a}_n$ and $\vec{b}_{2n}$ obtaining as output the reconstructed scaling vector $\vec{v}$ 
\end{algorithmic}

\end{algorithm}

\MVB{In Figure~\ref{fig_fsca_fwav} a comparison is shown between the choices
$f_{sca} = \sqrt{3}$, $f_{wav}=\sqrt{\frac 32}$ and $f_{sca} = f_{wav} = 1$.
Note that for $f_{sca} = \sqrt{3}$, $f_{wav}=\sqrt{\frac 32}$ the sum of
the squared norms is almost equal to $1$.
\begin{figure}[ht!]%
\begin{center}
\includegraphics[scale = 0.30]{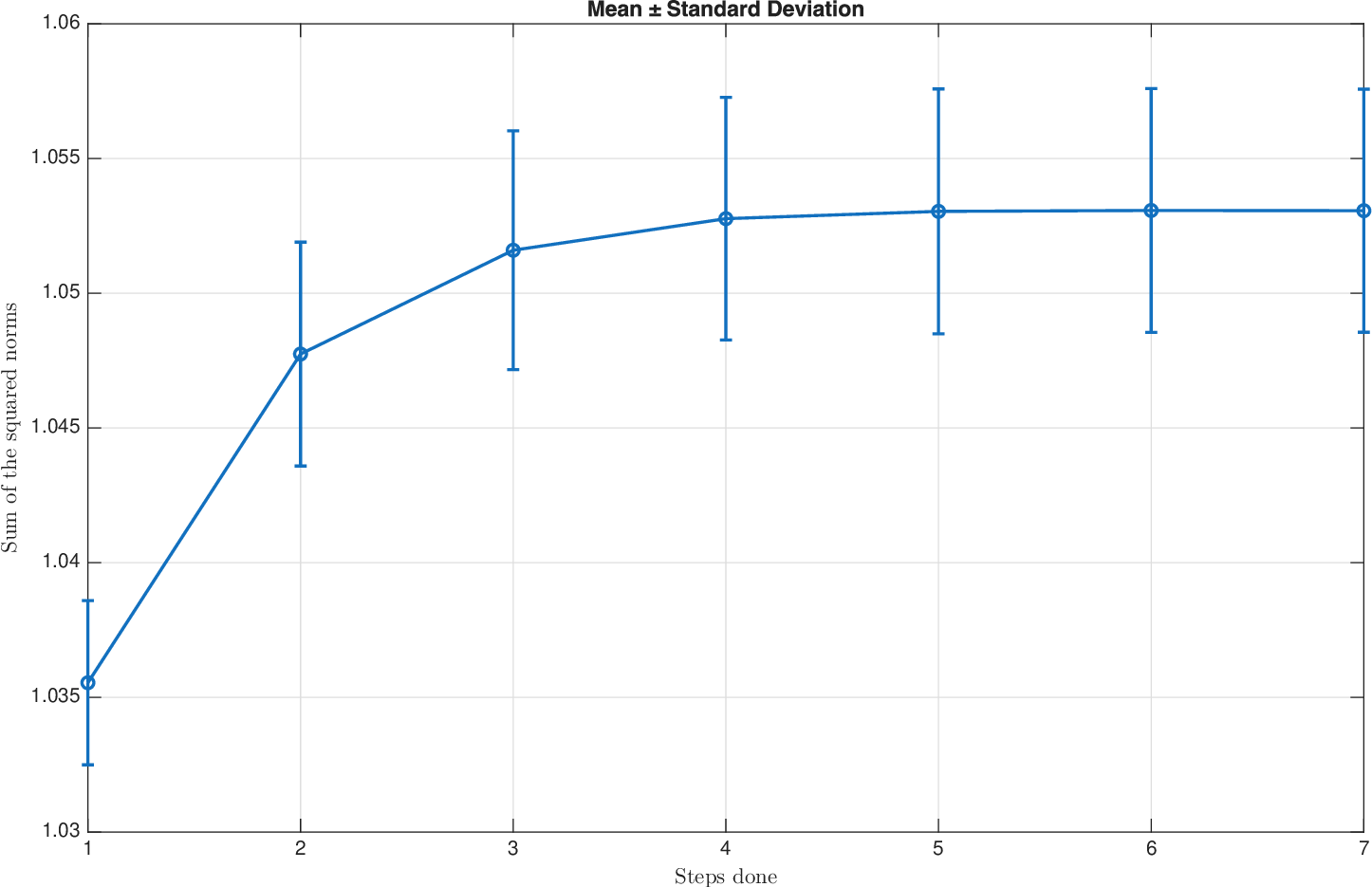}
\includegraphics[scale = 0.30]{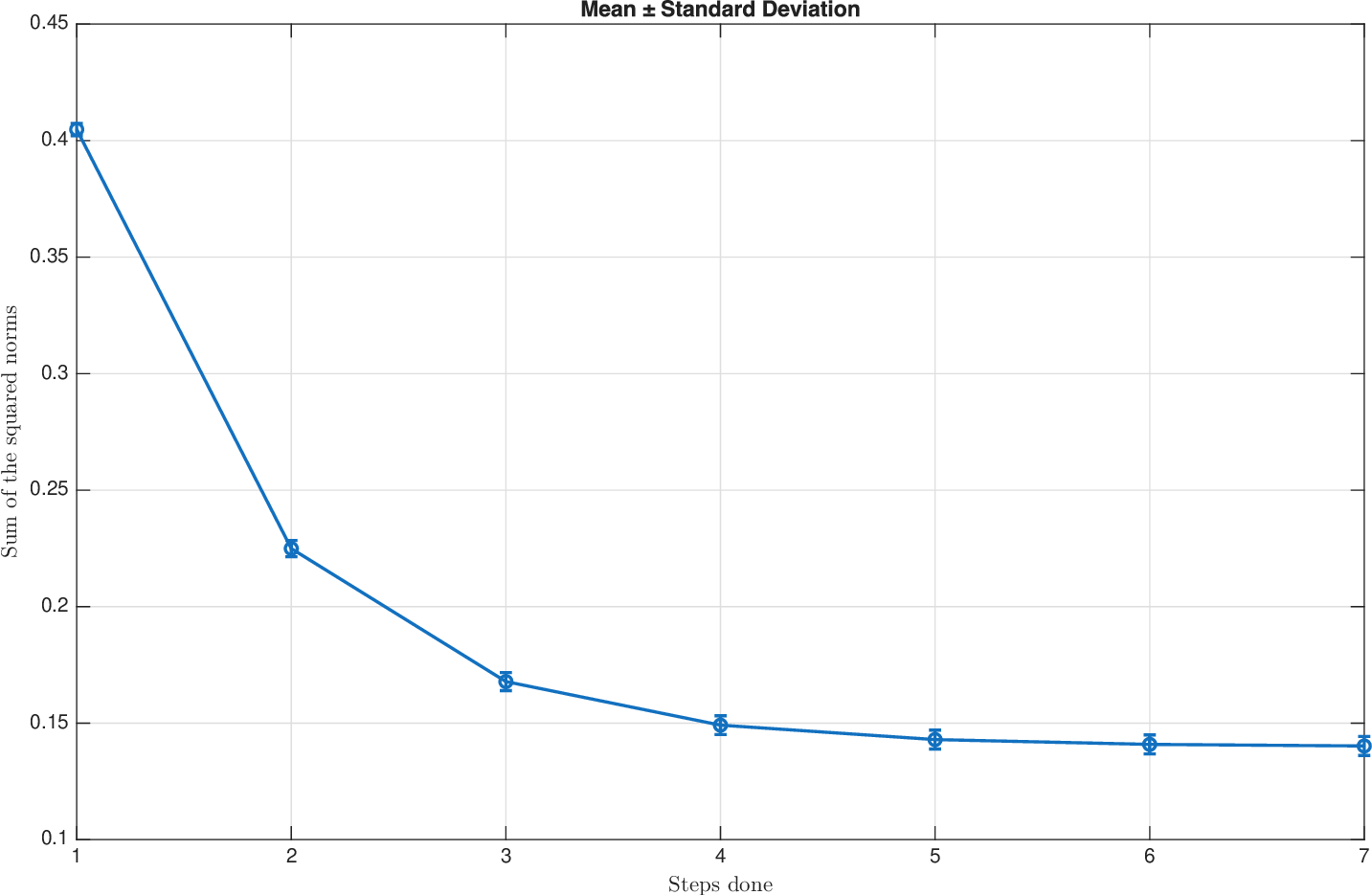}
\end{center}
\caption{ \MVB{The mean and standard deviation of the sum of the squared norms for $100$ samples for a signal of length $3^7$, $\theta = 0.5$ and $7$ steps 
done for $f_{sca} = \sqrt{3}$, $f_{wav}=\sqrt{\frac 32}$ (left) and for $f_{sca} = f_{wav} = 1$ (right).} \label{fig_fsca_fwav} }
\end{figure}
}

Let us conclude the section with some considerations concerning the extension of the previous arguments to the bivariate case that is  of interest in the case that the signal to process is an image. 

As is well--known, mathematically, any digital (gray-level) image at the resolution $N\times M$ is a matrix $I$ with size $N\times M$. According to the mathematical model already introduced in \cite{Image_VP2D},  we interpret its pixels $I_{h,k}$ as the values of an ideal function $f$ at the Chebyshev grid with the same size, namely
\[
I_{h,k}=f(x_h^N, x_k^M), \qquad \mbox{with}\quad
x_j^n = \cos\frac{(2j-1)\pi}{2n}, \quad j=1,\ldots,n.
\]
Hence, these pixels are the coefficients of the bivariate VP interpolation polynomial achieved via tensor product, i.e.
\begin{equation}\label{VP-2D}
V_{N,M}^{\nu,\mu}f(x,y)=\sum_{h=1}^N\sum_{k=1}^M I_{h,k}\ \phi_{N,h}^\nu(x)\phi_{M,k}^\mu(y),
\qquad  (x,y)\in [-1,1]^2
\end{equation}
where, as usual, we set $\nu=\lfloor\theta N
\rfloor$ and $\mu=\lfloor\theta M\rfloor$ with arbitrarily fixed $\theta\in ]0,1[$.\newline
Similarly to the univariate case, the VP interpolating polynomial $V_{N,M}^{\nu,\mu}f$ projects the unknown function $f$ onto the approximation space $\V_N^\nu\times \V_M^\mu$ and represents a near best polynomial approximation of $f$, as 
studied in \cite{OT-VPsquare}.

From the wavelets point of view, the $N\times M$ pixels $I_{h,k}$ constitute the scaling coefficients forming the scaling matrix $I$ to be decomposed and then reconstructed through a standard tensor product approach. \newline
Indeed, supposing for simplicity that both the dimensions $N$ and $M$ are divisible by three, e.g. setting
\[N=3N_1\qquad\mbox{and}\qquad M=3M_1
\] 
the approximation space $\V_N^\nu\times\V_M^\mu$ can be decomposed into the orthogonal sum of the lower-degree approximation space and three detail spaces corresponding to the horizontal, vertical, and diagonal directions, i.e.
\[
\V_N^\nu\times \V_M^\mu = \left(\V_{N_1}^\nu\times \V_{M_1}^\mu \right)\bigoplus \left( \W_{N_1}^\nu\times \V_{M}^\mu \right) \bigoplus \left(\V_{N}^\nu\times \W_{M_1}^\mu\right)\bigoplus \left(\W_{N_1}^\nu\times \W_{M_1}^\mu\right).
\]
This decomposition of spaces corresponds to a similar decomposition of the VP interpolating polynomial \eqref{VP-2D} as the sum of one lower-degree approximation and three detail polynomials. 
	In terms of coefficient matrices in the scaling and wavelet bases defined by tensor products, the initial scaling matrix $I$ can therefore be decomposed (and reconstructed) into a scaling matrix at the lower resolution $N_1\times M_1$ and three detail matrices of size $2N_1\times M$ (horizontal details), $N\times 2M_1$ (vertical details), and $2N_1\times 2M_1$ (diagonal details), respectively. 
    
Hence, the decomposition and reconstruction of a matrix $I$ can be carried out using the same one-dimensional algorithms (Algorithms 5 and 6) iteratively applied to the columns and rows of the image data

}
With the described algorithms and all the above observations in mind, we now proceed to test the VP basis through experiments on two typical data processing tasks: signal denoising and image compression. In the first case, we use univariate VP wavelets, in the second case we deal with bidimensional VP wavelets. 


	\begin{figure}[ht!]%
\begin{center}
\includegraphics[scale = 0.25]{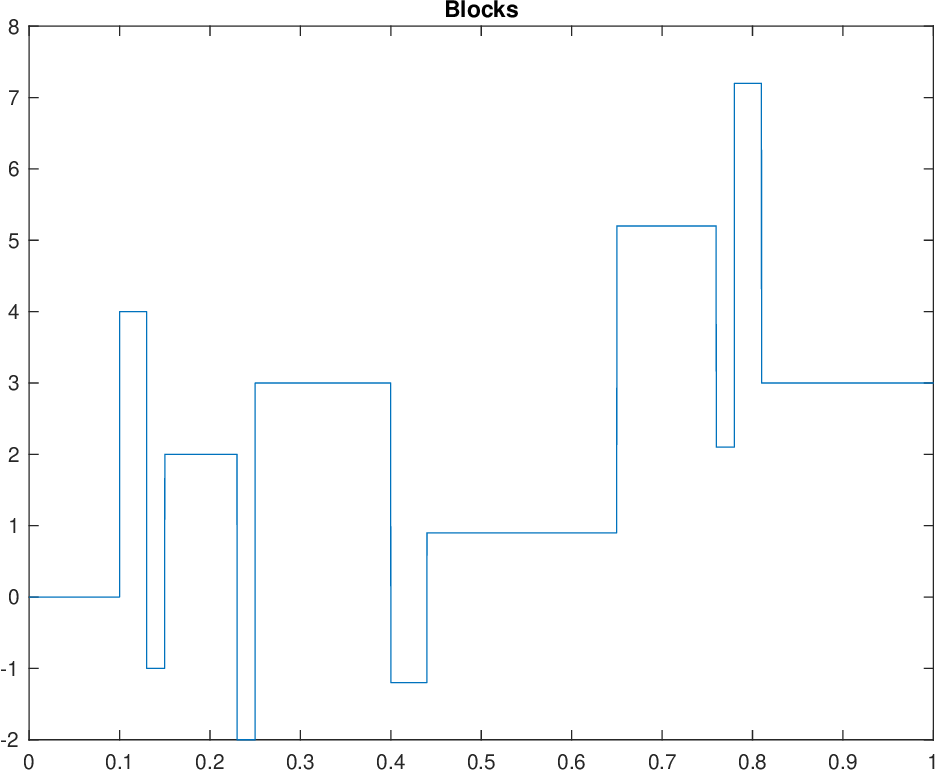}
\includegraphics[scale = 0.25]{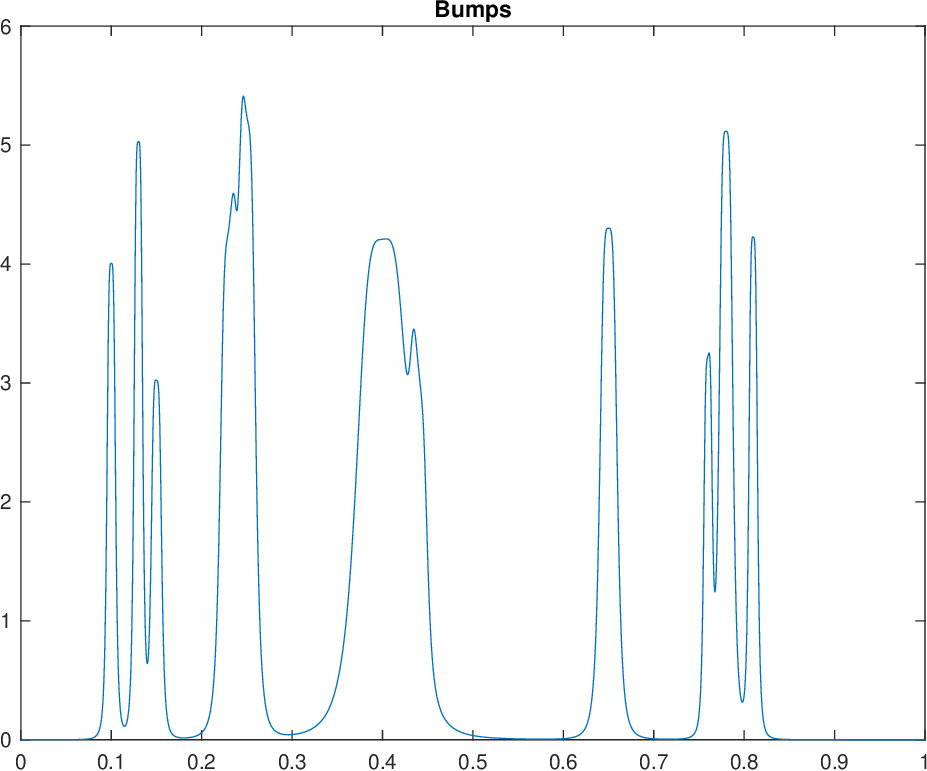}
\includegraphics[scale = 0.25]{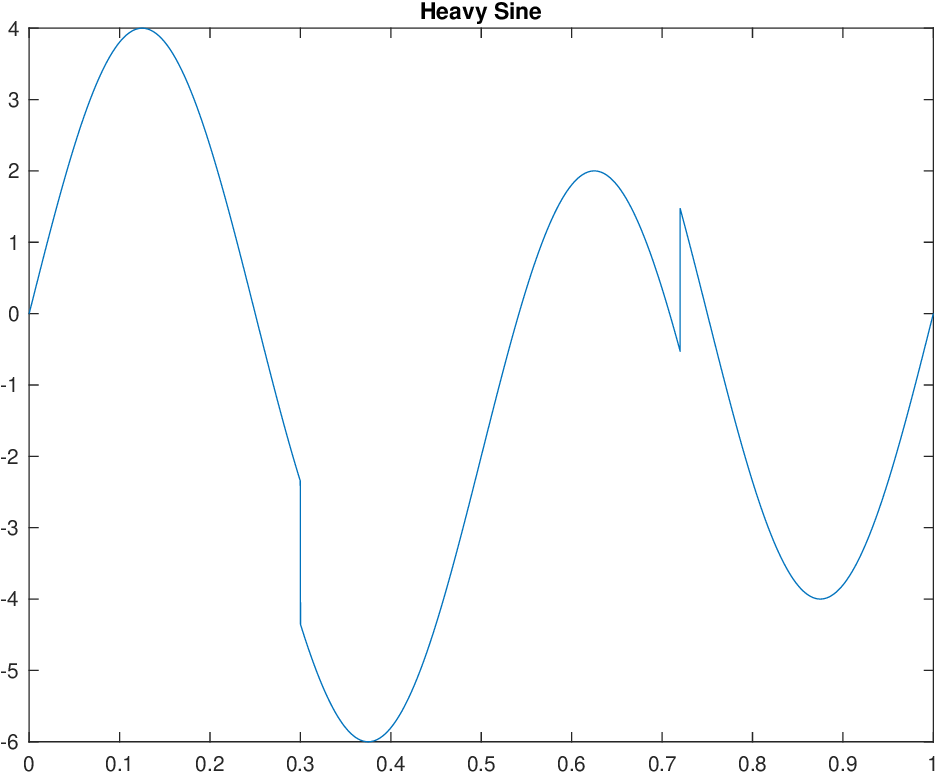}
\includegraphics[scale = 0.25]{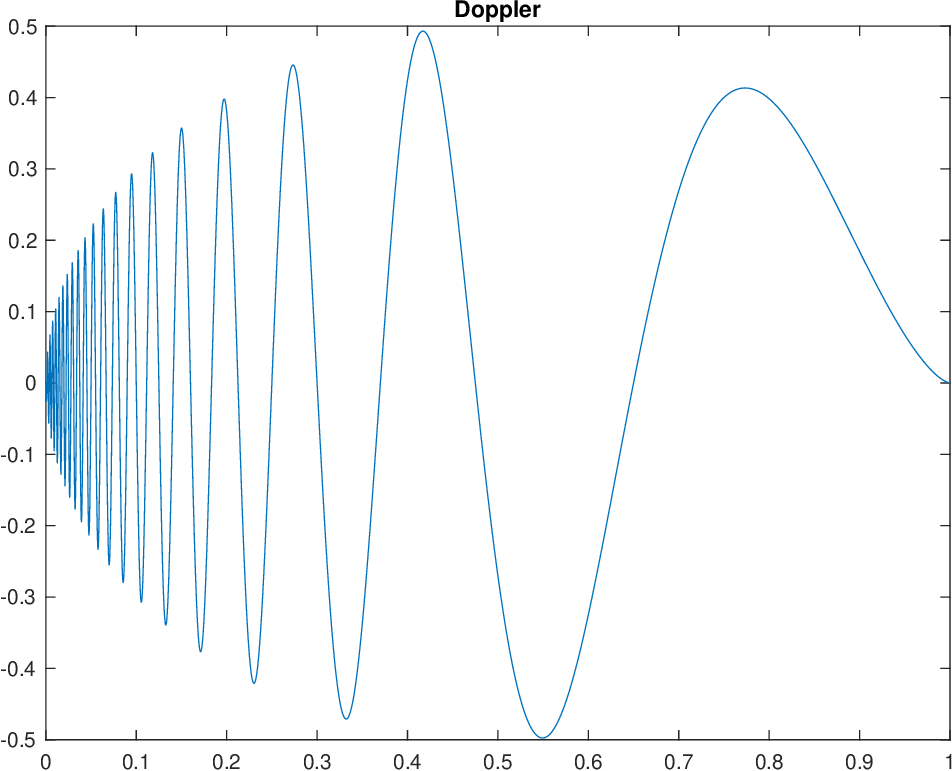}
\includegraphics[scale = 0.25]{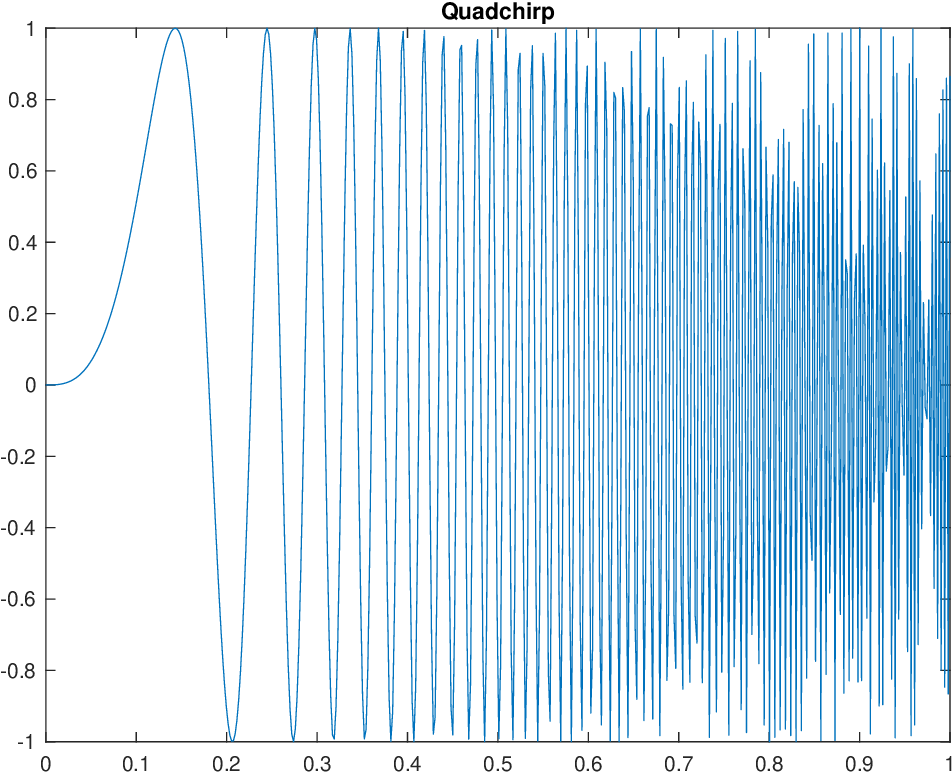}
\includegraphics[scale = 0.25]{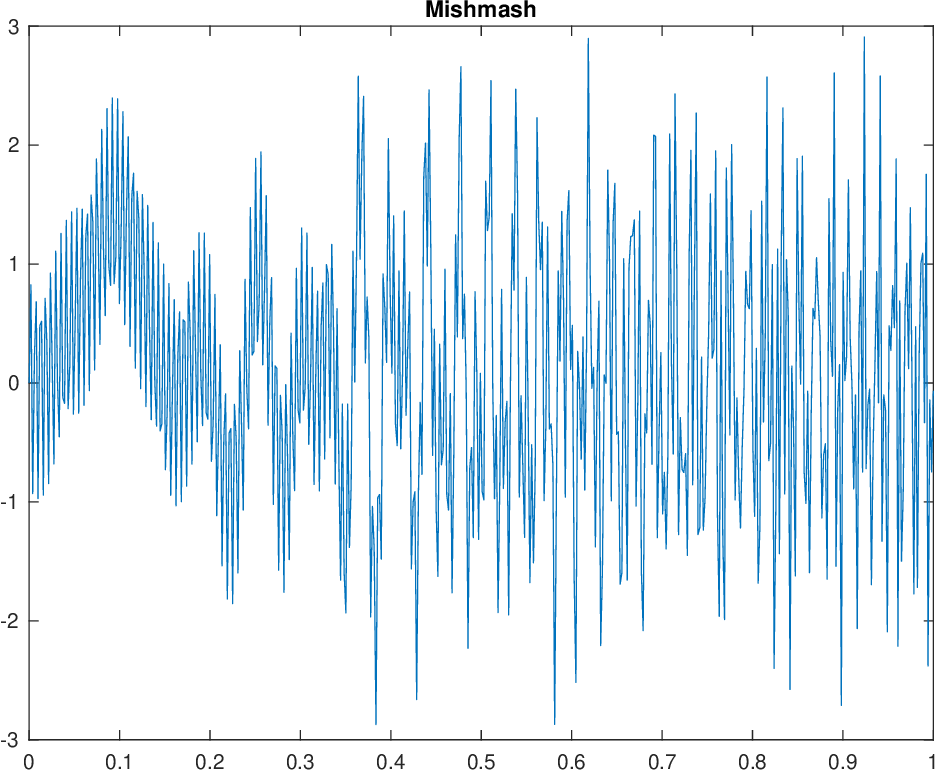}
\end{center}
\caption{ Six \MVB{test signals}: Blocks, Bumps, Heavy Sine, Doppler, Quadchirp, Mishmash.\label{fig_test_functions} }
\end{figure}

\begin{figure}[ht!]%
\begin{center}
\includegraphics[scale = 0.40]{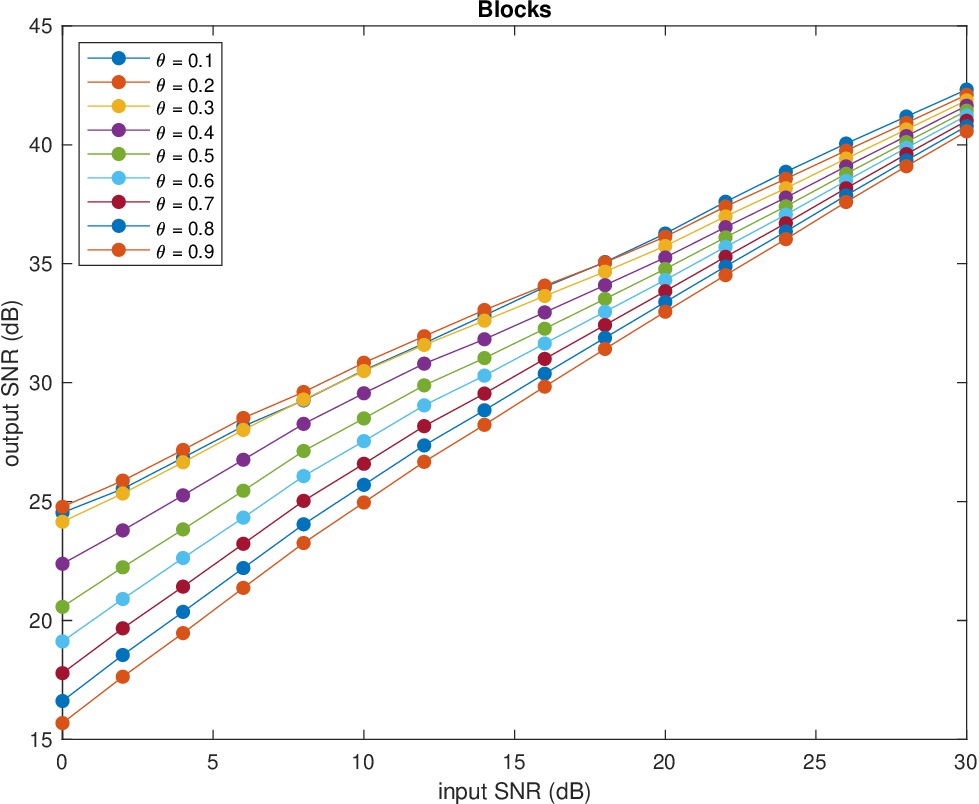}
\includegraphics[scale = 0.40]{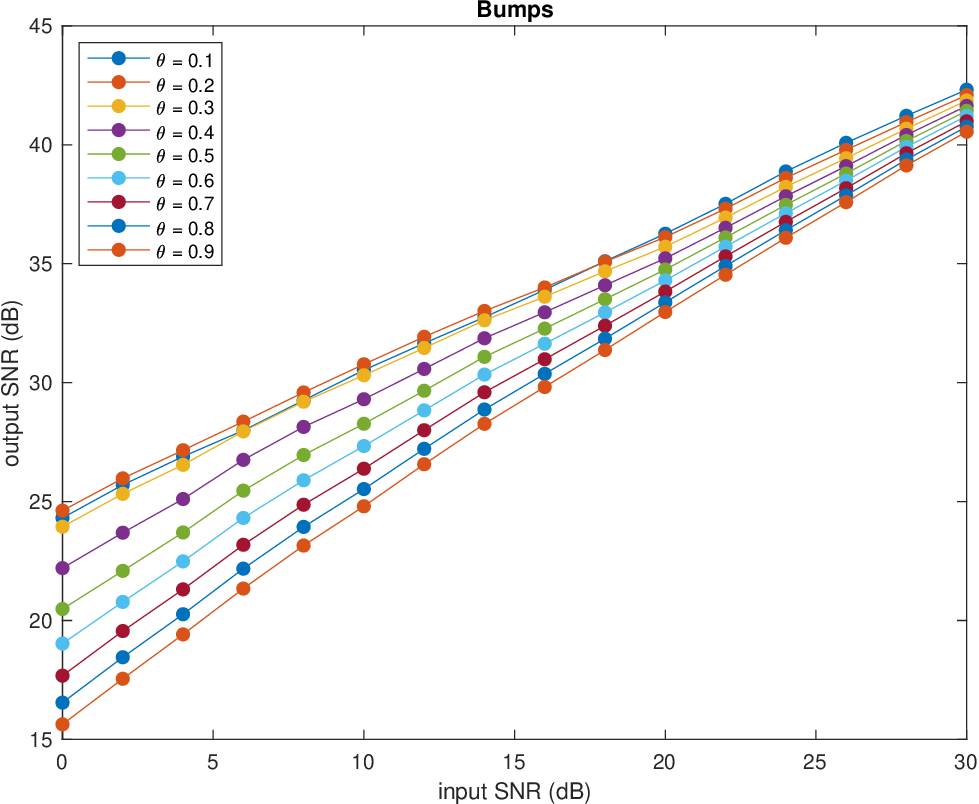}\\
\includegraphics[scale = 0.40]{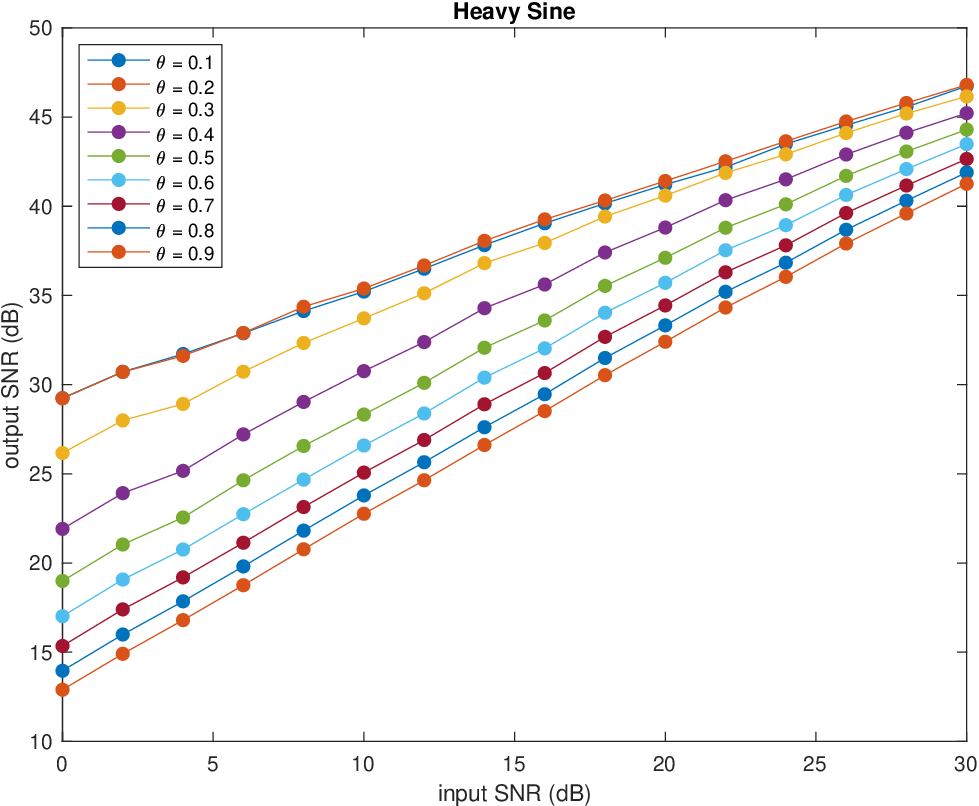}
\includegraphics[scale = 0.40]{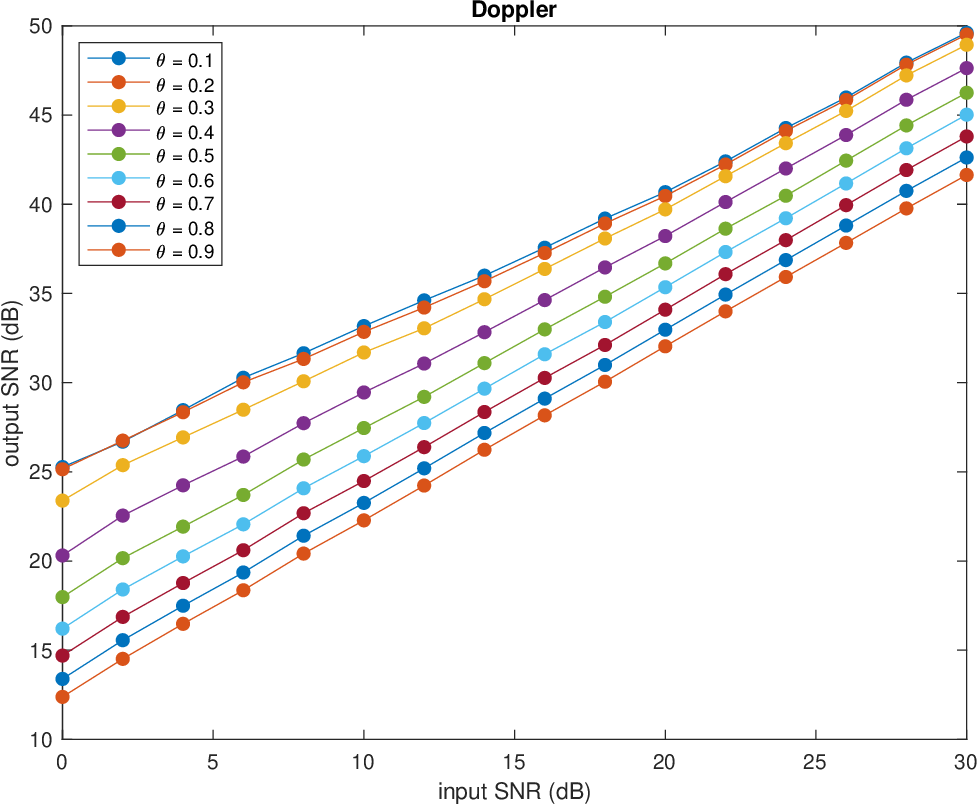}\\
\includegraphics[scale = 0.40]{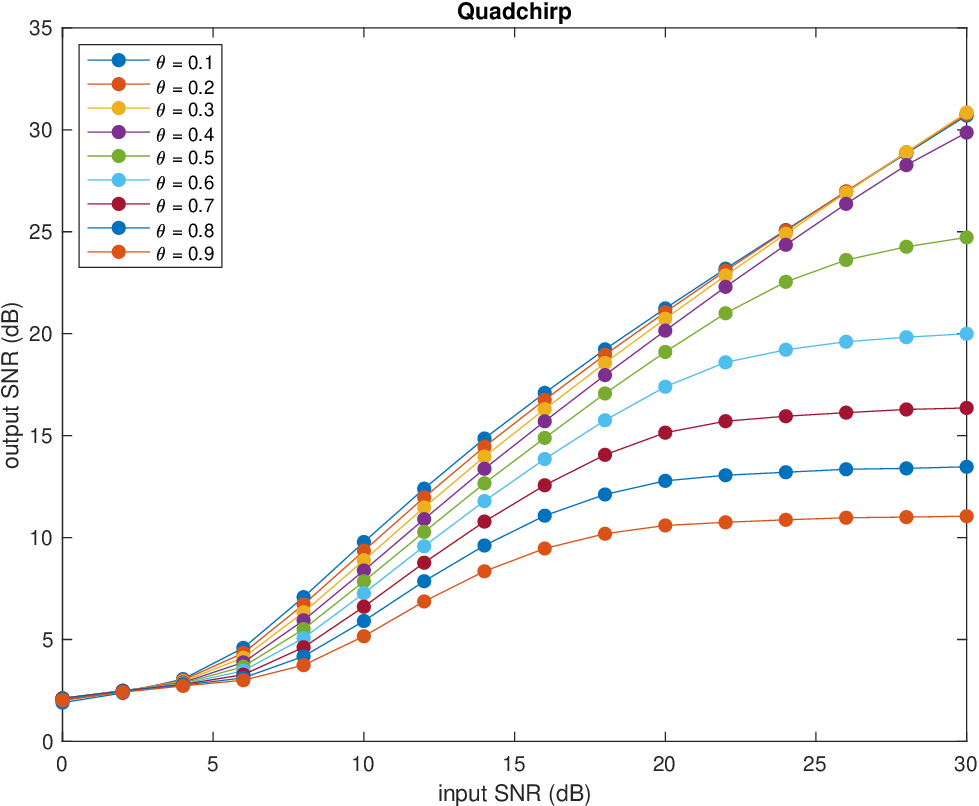}
\includegraphics[scale = 0.40]{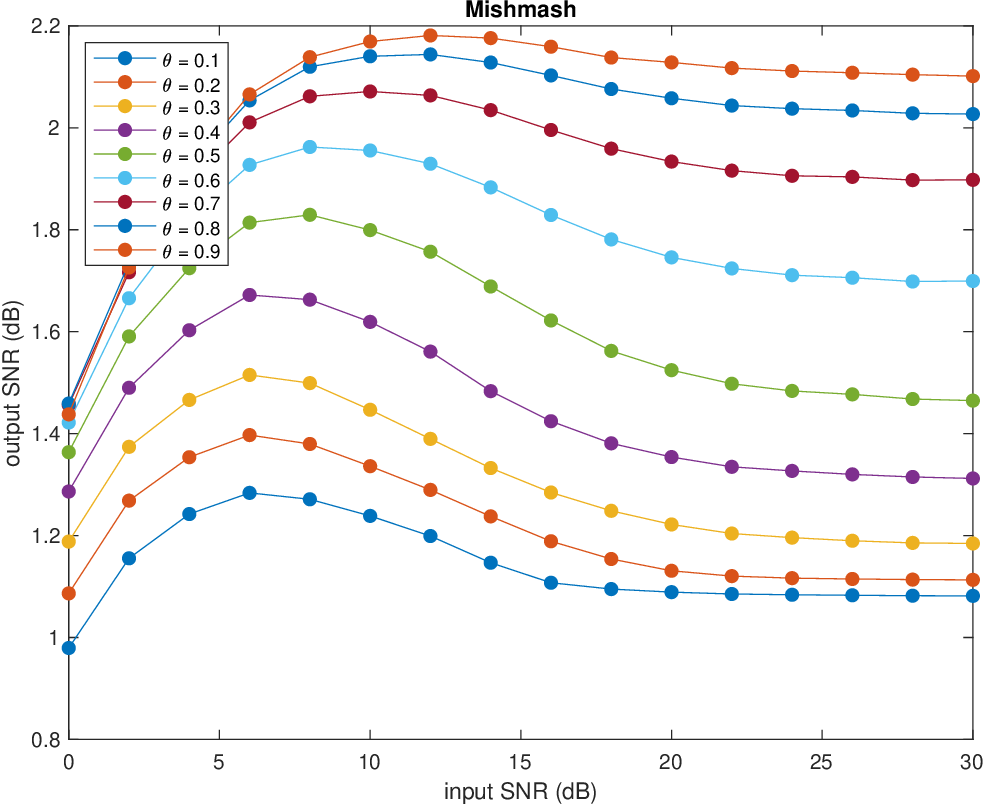}
\end{center}
\caption{ The results for different values of $\theta = 0.1, 0.2, \ldots, 0.9$ for the six \MVB{test signals}: Blocks, Bumps, Heavy Sine, Doppler, Quadchirp, Mishmash.\label{fig_theta} }
\end{figure}

\section{Signal denoising}
Wavelet transform-based  methods have proven to be highly effective in signal denoising tasks. 
The goal is to recover a signal $s= (s_k : k = 1, \dots, n)$ from noisy observations $y = (y_k : k = 1, \dots, n)$:
$$
y_k = s_k + \eta_k, \quad k = 1, \dots, n,
$$
where $\eta_k \sim \mathcal{N}(0, \sigma^2)$ are independent and identically distributed (i.i.d.) Gaussian random variables with mean 0 and variance $\sigma^2$.


Several different algorithms for wavelet-based signal denoising have been proposed in the literature, starting with the pioneering works by Donoho and Johnstone \cite{donoho1994ideal,donoho1995denoising}. Each method applies a different thresholding strategy to separate the signal from the noise.

In this context, we aimed at testing such methods  when used in conjunction with VP wavelets, assessing their performance in comparison with standard wavelet filters. In our preliminary experiments, it was found that the Bayesian method \cite{chipman1997adaptive} with soft thresholding produced the best results. 

\MVB{Using the factors $\sqrt{3}$ and $\sqrt{3/2}$ as proposed at the end of the previous section already improved  the results but we wanted to refine these normalization factors.
Therefore,
we took 1000 Gaussian random signals with mean $0$ and variance $1$.
These signals were decomposed and the standard deviation for each of the detail
components and for the approximation component was measured.
In case of orthogonal wavelets, the averages of these standard deviations
would converge to one if the number of samples increases.
For the non-orthogonal wavelets used in this paper, the magnitude of these averages turned out to be around $1$.
We used these averages as normalization factors.}


To determine how the parameter $\theta$ influences the reconstruction results,  we applied the Bayes method to the six \MVB{test signals} of the
wavelet toolbox of Matlab\footnote{We used Matlab version R2024B.
Matlab is a registered trademark of the MathWorks.}: Blocks, Bumps, Heavy Sine, Doppler, Quadchirp and Mishmash,
plotted in Fig.~\ref{fig_test_functions}.

For each of these \MVB{test signals} we executed the Bayes method with soft thresholding 
for parameter values
$\theta$ ranging from $0.1$ to $0.9$. The length of the input signals was $2^{17}$.
\MVB{The depth of the iterations was $8$.}
For each value of $\theta$ $10$ samples of the input signal were considered
and the average was taken \MVB{as well as the standard deviaton with respect
to the average.}

\begin{figure}[!htb]%
\begin{center}
\includegraphics[scale = 0.40]{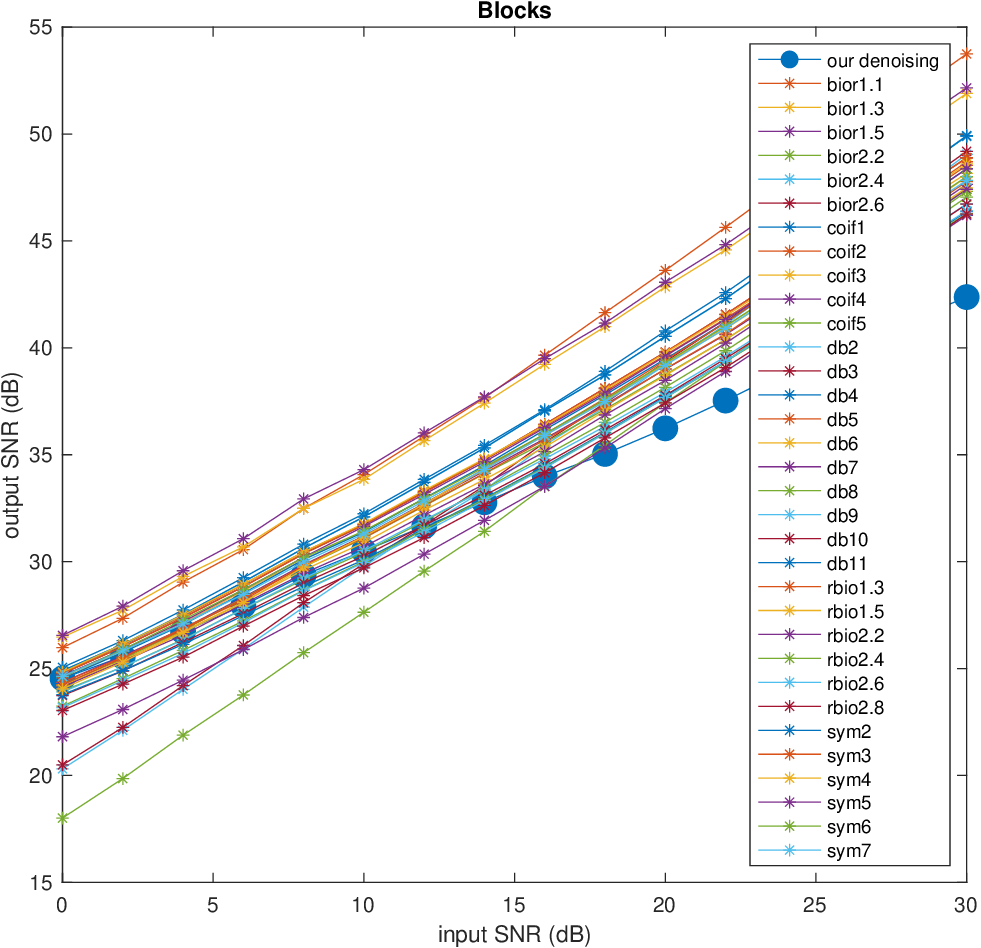}
\includegraphics[scale = 0.40]{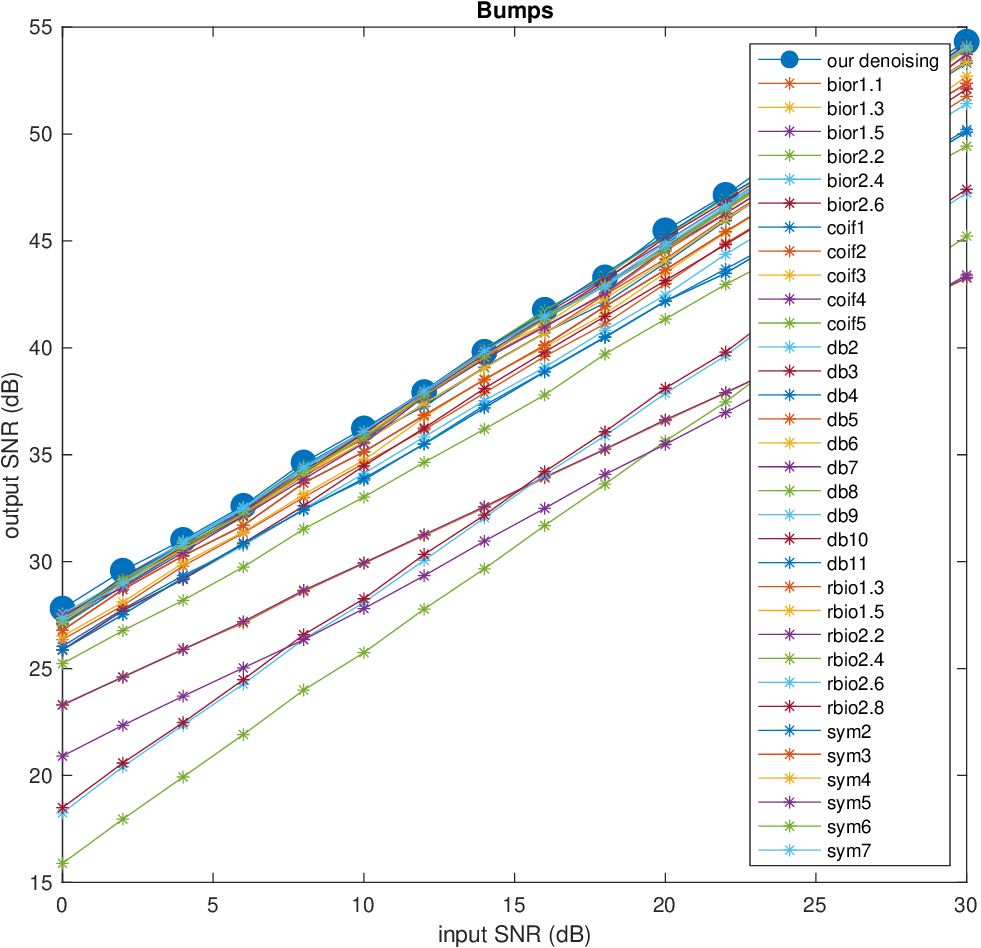}\\
\includegraphics[scale = 0.40]{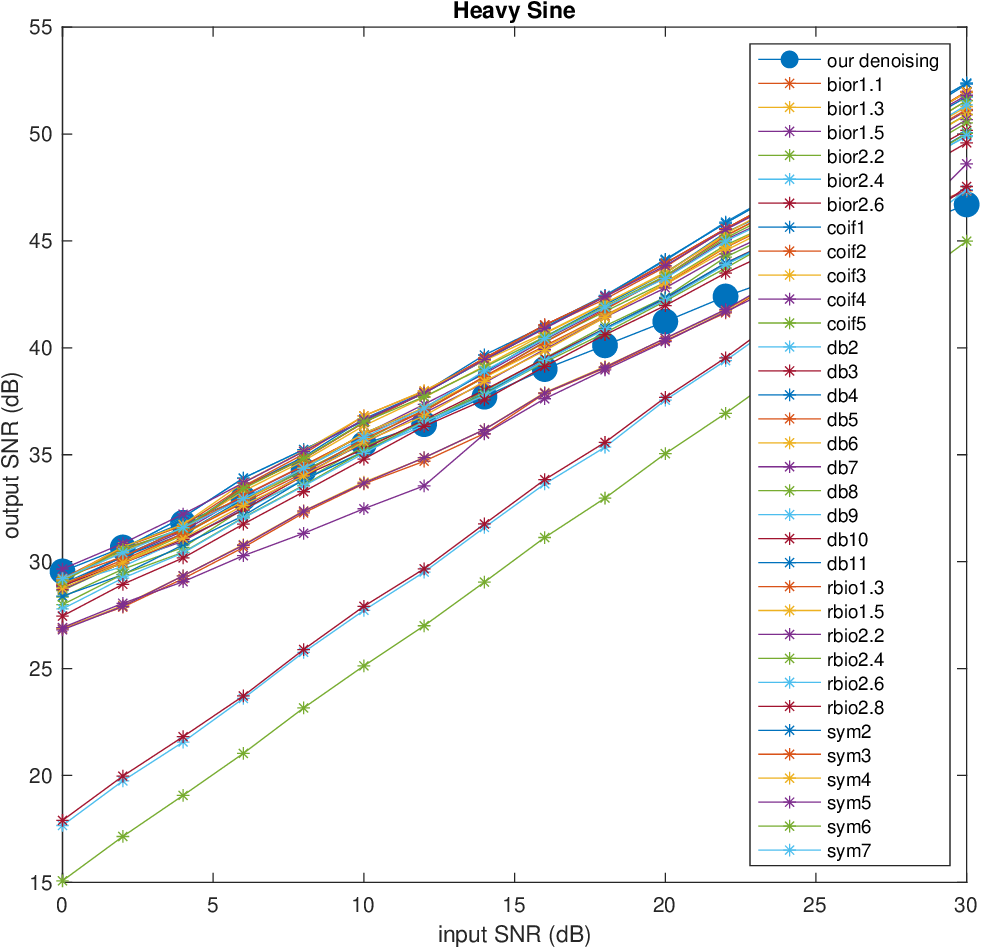}
\includegraphics[scale = 0.40]{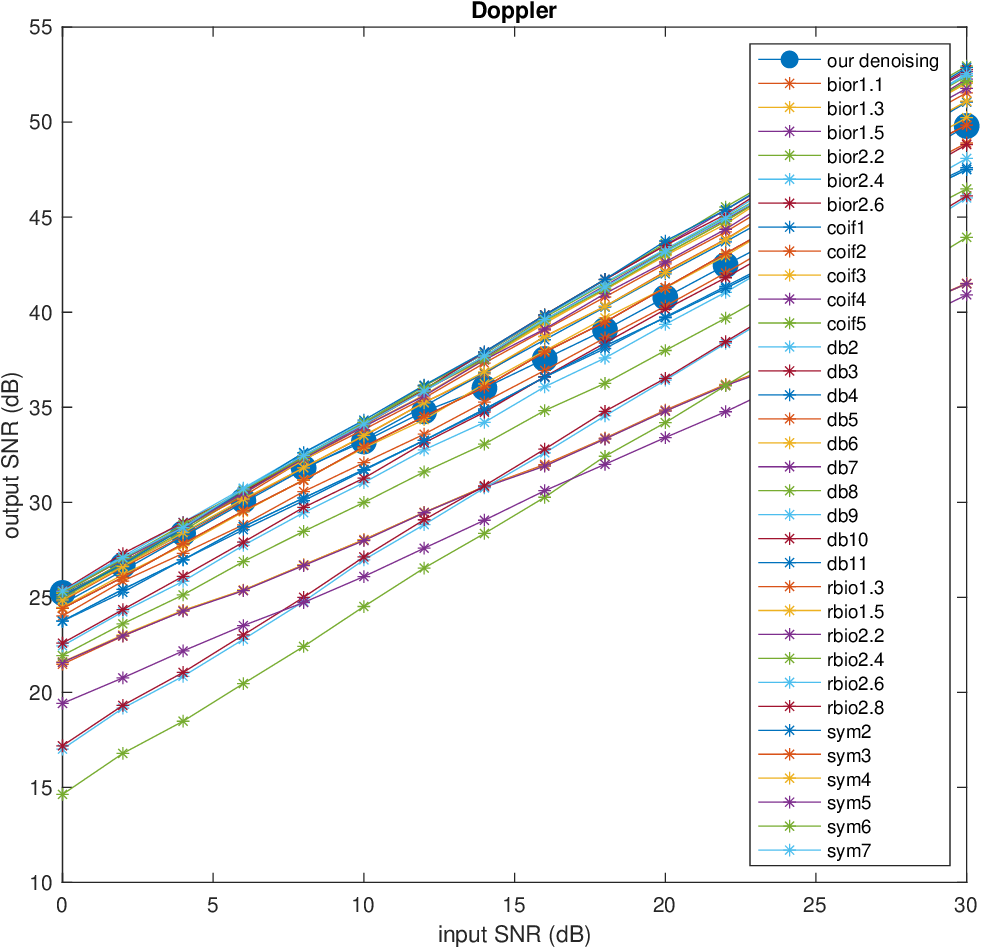}\\
\includegraphics[scale = 0.40]{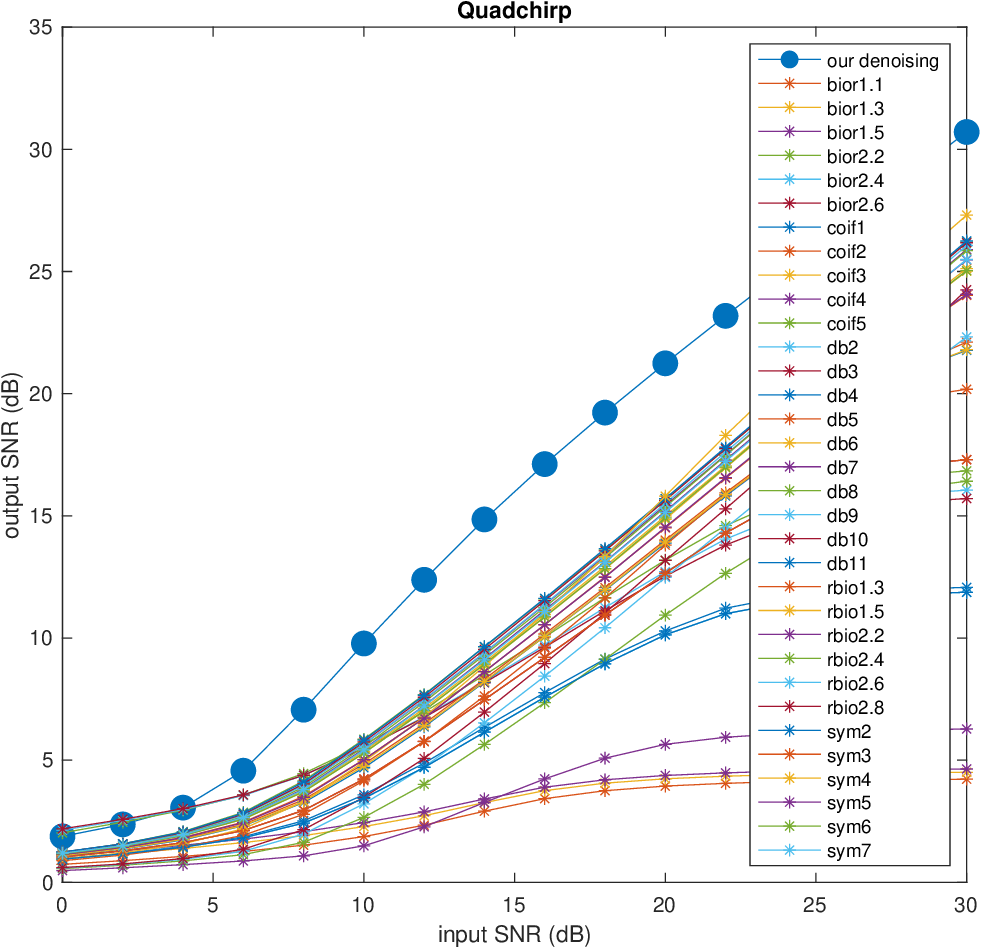}
\includegraphics[scale = 0.40]{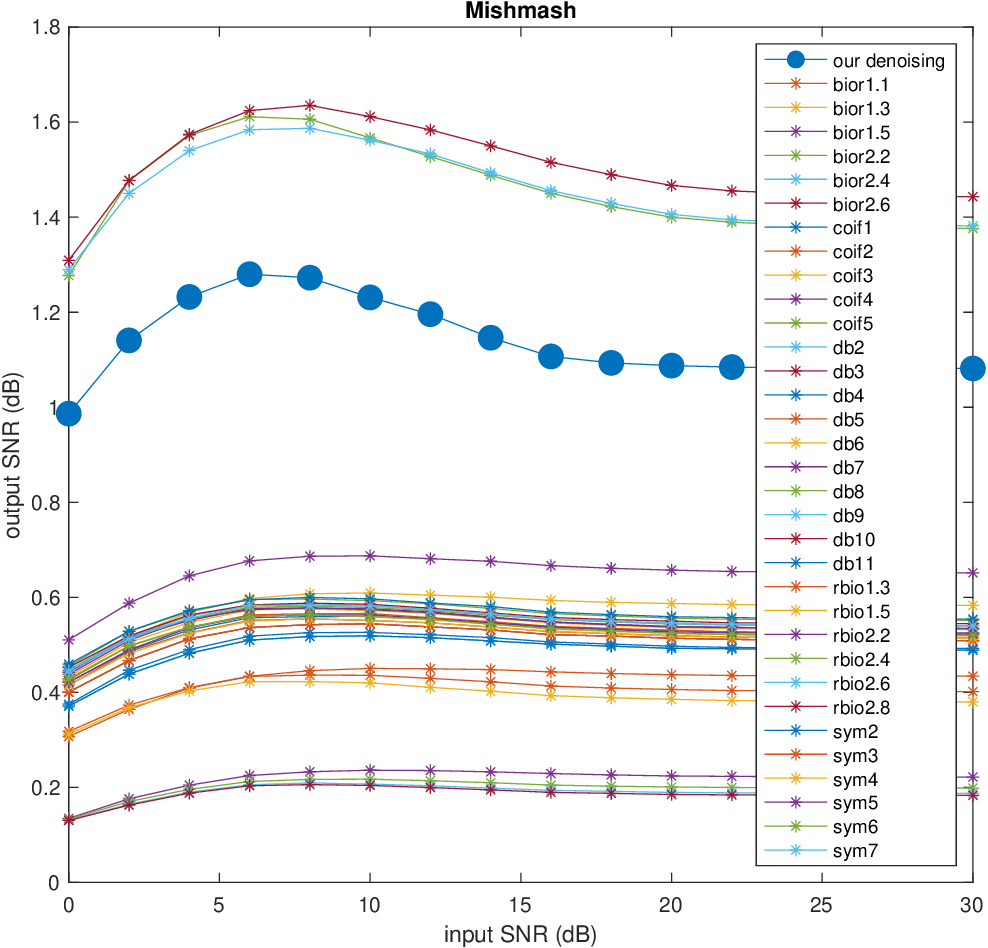}
\end{center}
\caption{ Comparison of our denoising method with the ones of ``wdenoise'' of Matlab
for the six \MVB{test signals}: Blocks, Bumps, Heavy Sine, Doppler, Quadchirp, Mishmash.\label{fig_denoise} }
\end{figure}

In Fig.~\ref{fig_theta}, the result is plotted for the six \MVB{test signals}.
This figure shows the signal-to-noise ratio (SNR) of the input signal with respect 
to the SNR of the output signal. We recall that the SNR of the two signals $s_1$ and $s_2$ is defined as
\[
\mbox{SNR}(s_1,s_2) = 10 \log_{10} \left( \frac{\sigma_{signal}^2}{\sigma_{noise}^2}
\right)
\]
with $\sigma_{signal}^2$ the variance of the signal $s_1$ and $\sigma_{noise}^2$
the variance of the noise signal, i.e., the difference between the signals
$s_1$ and $s_2$.
The input SNR is defined as the SNR of the original signal and the noisy signal while
the output SNR is defined as the SNR of the original signal and the denoised signal.

\begin{figure}
\begin{center}
\includegraphics[width=10cm]{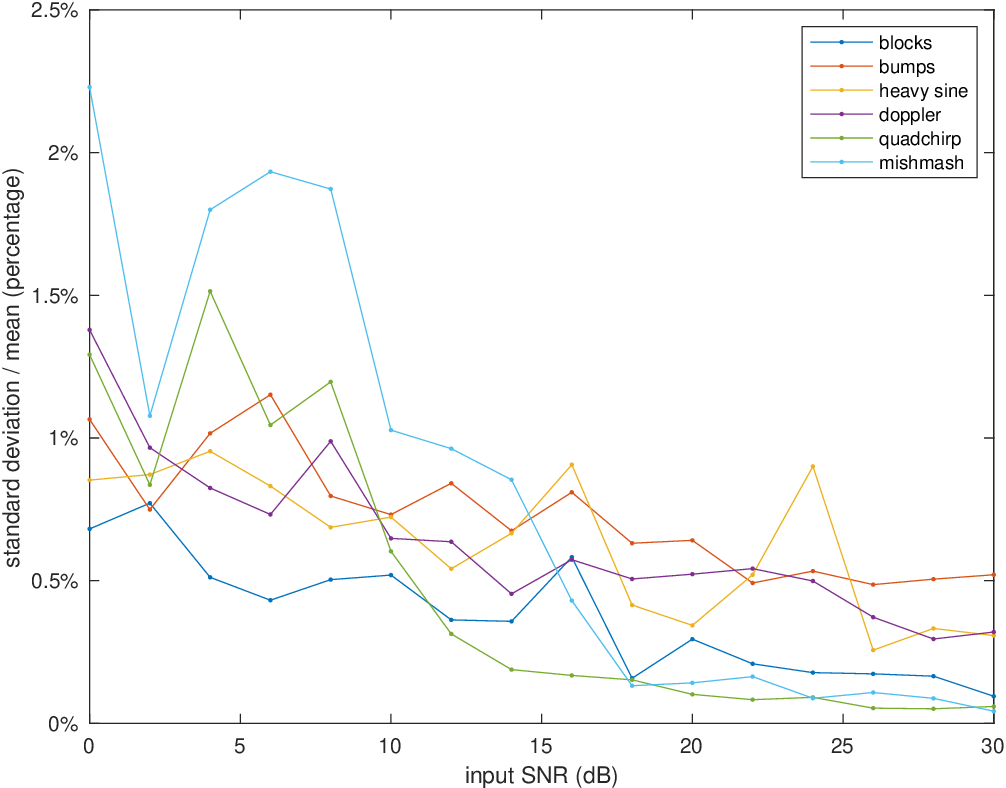}
\end{center}
\caption{\MVB{The standard deviation divided by the average in function of input SNR (dB) for each of the test
function for $\theta = 0.1$\label{fig:std} }}
\end{figure}

From the plots it is clear that the value $\theta = 0.1$ is best or near-best for the 6 \MVB{test signals}. Hence, in the sequel, we fixed the value of $\theta$ to be equal
to $0.1$.
\MVB{In Figure~\ref{fig:std} the standard deviation divided by the average is plotted for $\theta = 0.1$ in function of the input SNR (dB)
and each of the \MVB{test signals}.}

In the next experiment, we compared our denoising method with the one of Matlab ``wdenoise'' (of the wavelet toolbox) for all choices of the wavelets that are
available using this denoising function:
bior1.1, bior1.3, bior1.5, bior2.2, bior2.4, bior2.6,
        coif1, coif2, coif3, coif4, coif5,
        db2, db3, db4, db5, db6, db7, db8, db9, db10, db11,
        rbio1.3, rbio1.5, rbio2.2, rbio2.4, rbio2.6, rbio2.8,
        sym2, sym3, sym4, sym5, sym6, sym7.
Applying ``wdenoise'', we used the default denoising method which is ``Bayes'' and soft thresholding. Also here we took the average of $10$ samples.
Figure~\ref{fig_denoise} shows the results for the 6 \MVB{test signals}.

For ``Blocks'', our method is not performing well compared to the others.
In this case the best three wavelets are bior1.1, bior1.5 and bior1.3.
For ``Bumps'' our method is the best one for all values of the input SNR.
For ``Heavy Sine'' our method is the best for small values of the input SNR.
For ``Doppler'' also our method is one of the best for small values of the input SNR.
For ``Quadchirp'' our method is the best for almost all values of the input SNR.
For small values of the input SNR, bior2.6 is slightly better.
For ``Mishmash'' our method performs well with respect to the other wavelets.
Only 3 wavelets perform better: bior2.6, bior2.2 and bior2.4.

Overall, the performance of our basis varies depending on the signal type, nevertheless the results highlight the robustness of our approach  particularly in low-SNR conditions.

\begin{figure}
\begin{center}
\includegraphics[width=5cm]{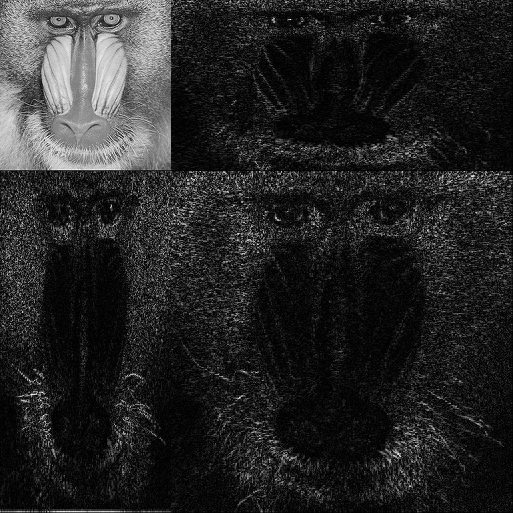}
\includegraphics[width=5cm]{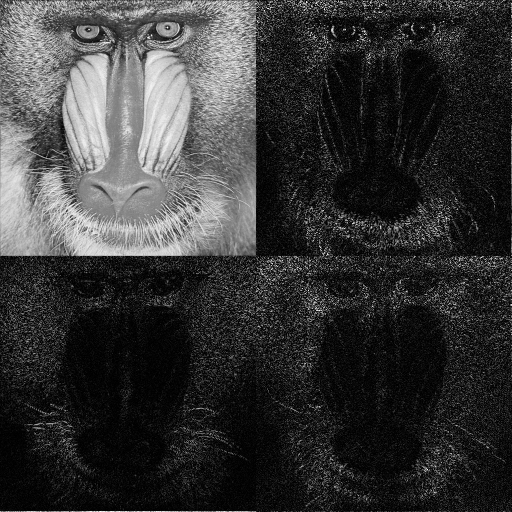}
\end{center}
\caption{1-level decomposition of the Baboon image with the tensor product VP wavelets (left) and the tensor product Daubechies wavelets (right) \label{fig:1stepdec} }
\end{figure}

	\section{Image compression}
Image compression is a key application of wavelet transforms, aimed at reducing storage and transmission costs while maintaining the perceptual quality of the data  (see, e.g., \cite{alarcon,BOIX, Bruni2020,chang,CRS2025,Usevitch, Ma}, among many others). In this section, we illustrate the performance of VP wavelets for image compression, comparing them with traditional separable wavelets. 
It is important to clarify that we use here the term "compression" somewhat loosely. Our approach involves only thresholding, without applying the full compression pipeline, which would typically include quantization and encoding. The goal of our study is in fact to assess the sparsifying ability - or compressibility - of the basis rather than to implement a complete compression scheme.

\begin{figure}[t]
\begin{center}
\includegraphics[width=7.2cm]{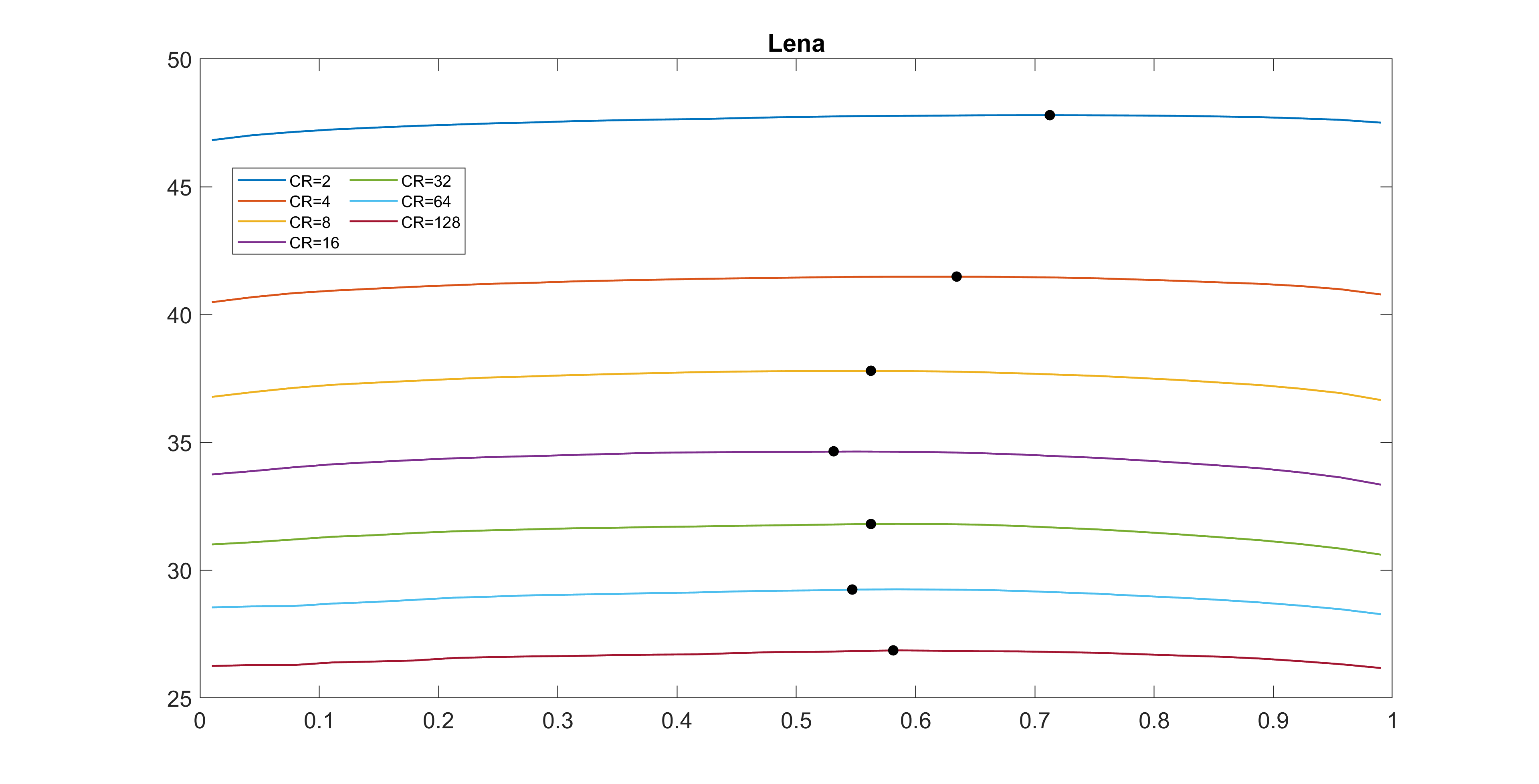}
\includegraphics[width=7.2cm]{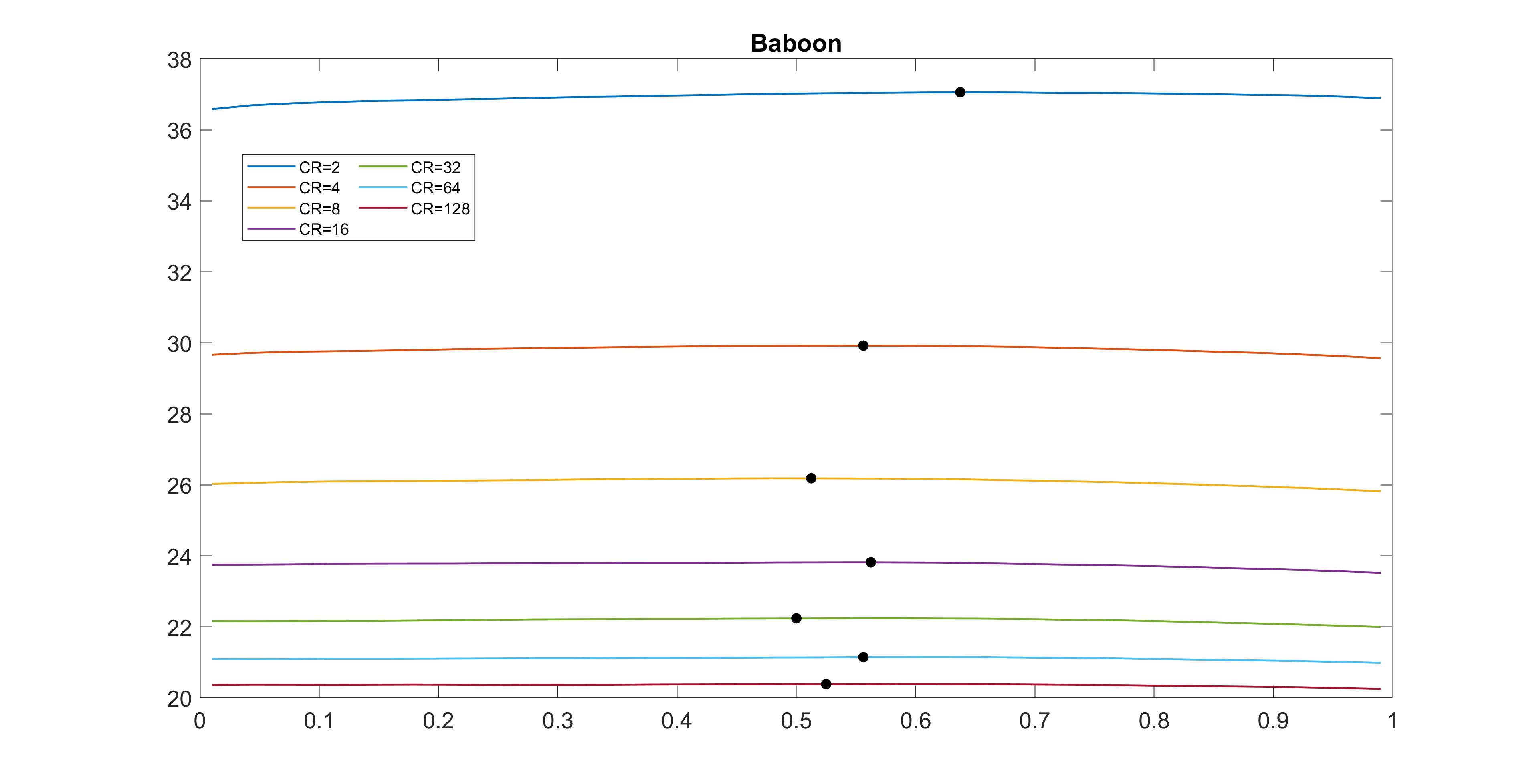}
\includegraphics[width=7.2cm]{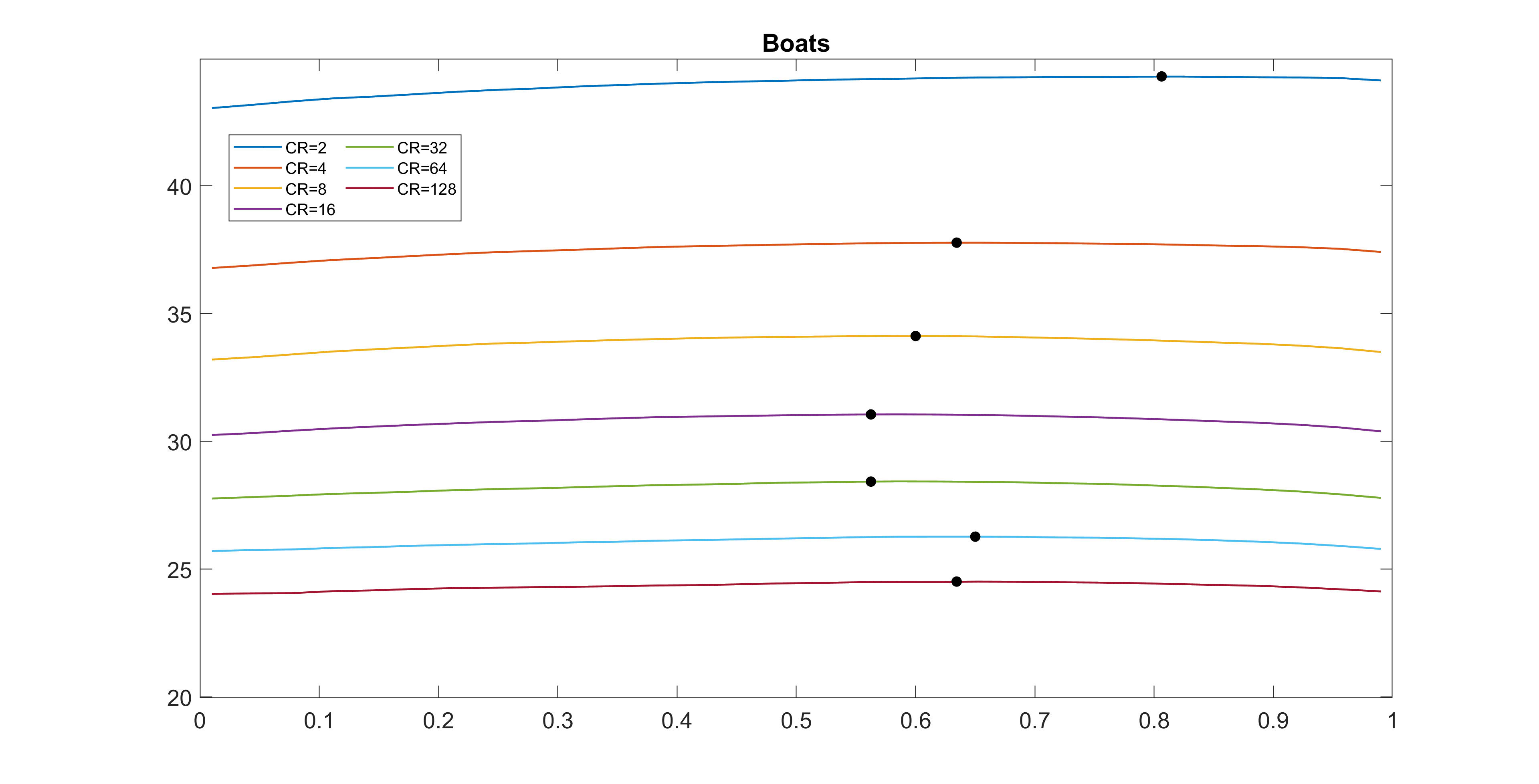}
\includegraphics[width=7.2cm]{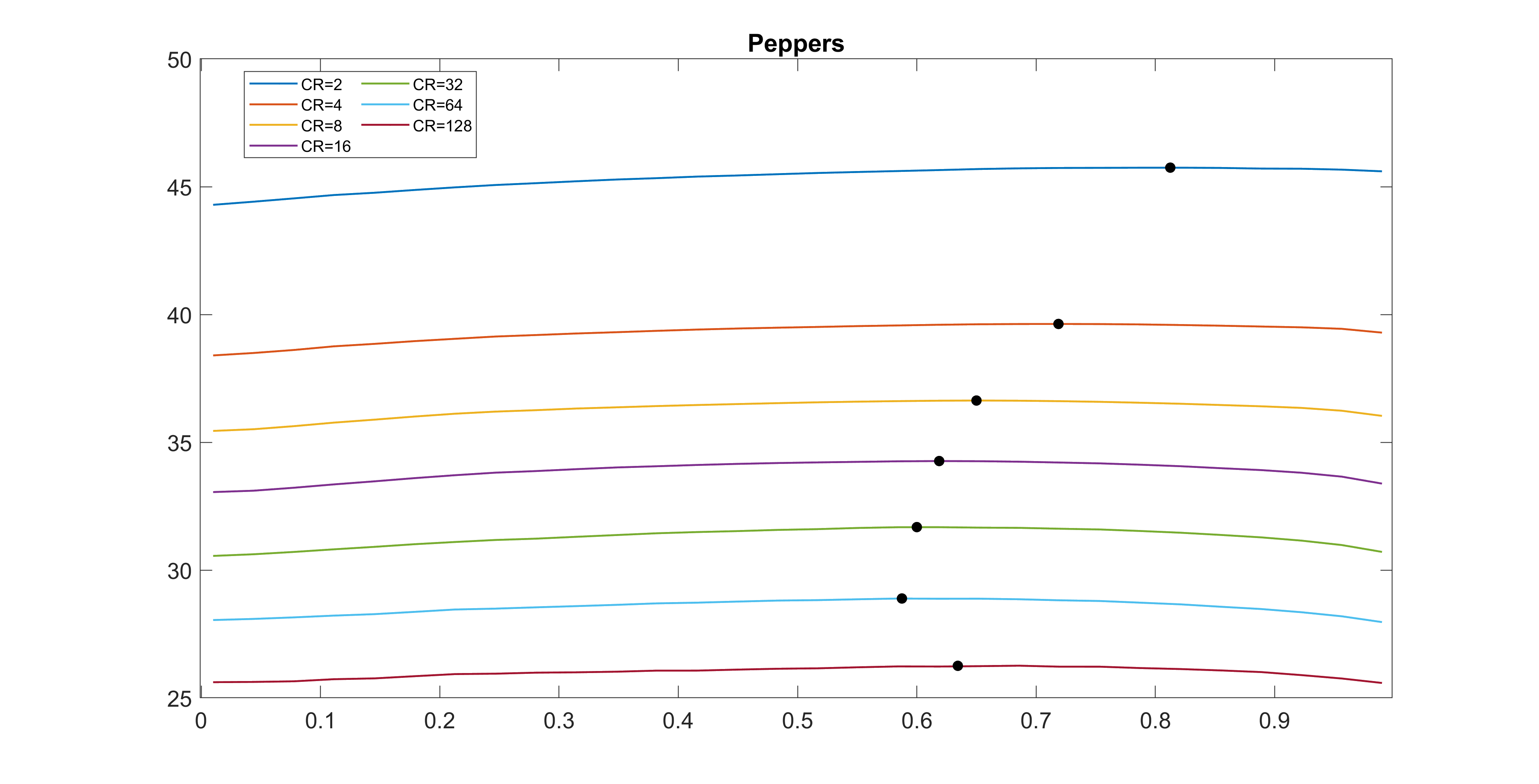}
\end{center}
\caption{Variation of PSNR as a function of the free parameter $\theta$ for different images and compression ratios. The optimal $\theta$ for each curve is indicated by a black bullet. \label{fig:PSNR_Theta} }
\end{figure}

We present several tests, focusing on decomposition efficiency, the trade-offs in compression ratio (CR), 
defined as the ratio between the size of the original image
and the size of the compressed image, and the impact of the parameter $\theta$.



We first wish to  illustrate the effect of a one-level decomposition using the VP basis. Fig. \ref{fig:1stepdec} presents the results of decomposing the classical Baboon image, comparing it to a standard wavelet tensor product decomposition (e.g., Daubechies 2). After just one level of decomposition,  the approximation component retains only 1/3 of the original data size, with the remaining portion consisting of detail coefficients. This reduction means that fewer decomposition steps are required compared to the standard dyadic transform, for achieving a similar compact image representation, leading to significant computational savings.

To evaluate the quality of the VP basis for image compression, we conducted standard experiments involving the following steps: selecting a compression ratio (CR), or, equivalently, the percentage P=100/CR of elements to be retained in the compressed image; decomposing the image up to a certain level; zeroing wavelet coefficients below a threshold determined by the CR; and reconstructing the image. 

The quality of the reconstruction is evaluated using both PSNR (Peak Signal-to-Noise Ratio) and SSIM (Structural Similarity Index Measure) metrics, and the results are compared to those achieved by using 
either the Daubechies filters with 4 coefficients (Db2) or the biorthogonal 3.5 (bior3.5) wavelets.

For standard wavelet decomposition and reconstruction, we used the MATLAB functions \texttt{wavedec2} and \texttt{waverec2}, respectively. The PSNR and SSIM quality metrics were computed using the MATLAB functions \texttt{psnr} and \texttt{ssim}.

The first tests were performed on four standard benchmark images: Lena, Baboon, Boats, and Peppers (each of size 512 $\times$ 512).  

The dependency of the VP basis on the free parameter $\theta$ allows for its optimal selection in image reconstruction. To illustrate the effect of $\theta$, Fig.~\ref{fig:PSNR_Theta} shows how the PSNR varies as a function of this parameter. This dependency is not uniform across all cases; rather, as observed, it is strongly influenced by the image data and the chosen compression ratio, which ranges as $2^k$ for $k = 1, \dots, 7$.

\begin{figure}[t]
 \begin{center}
 \includegraphics[width=6cm]{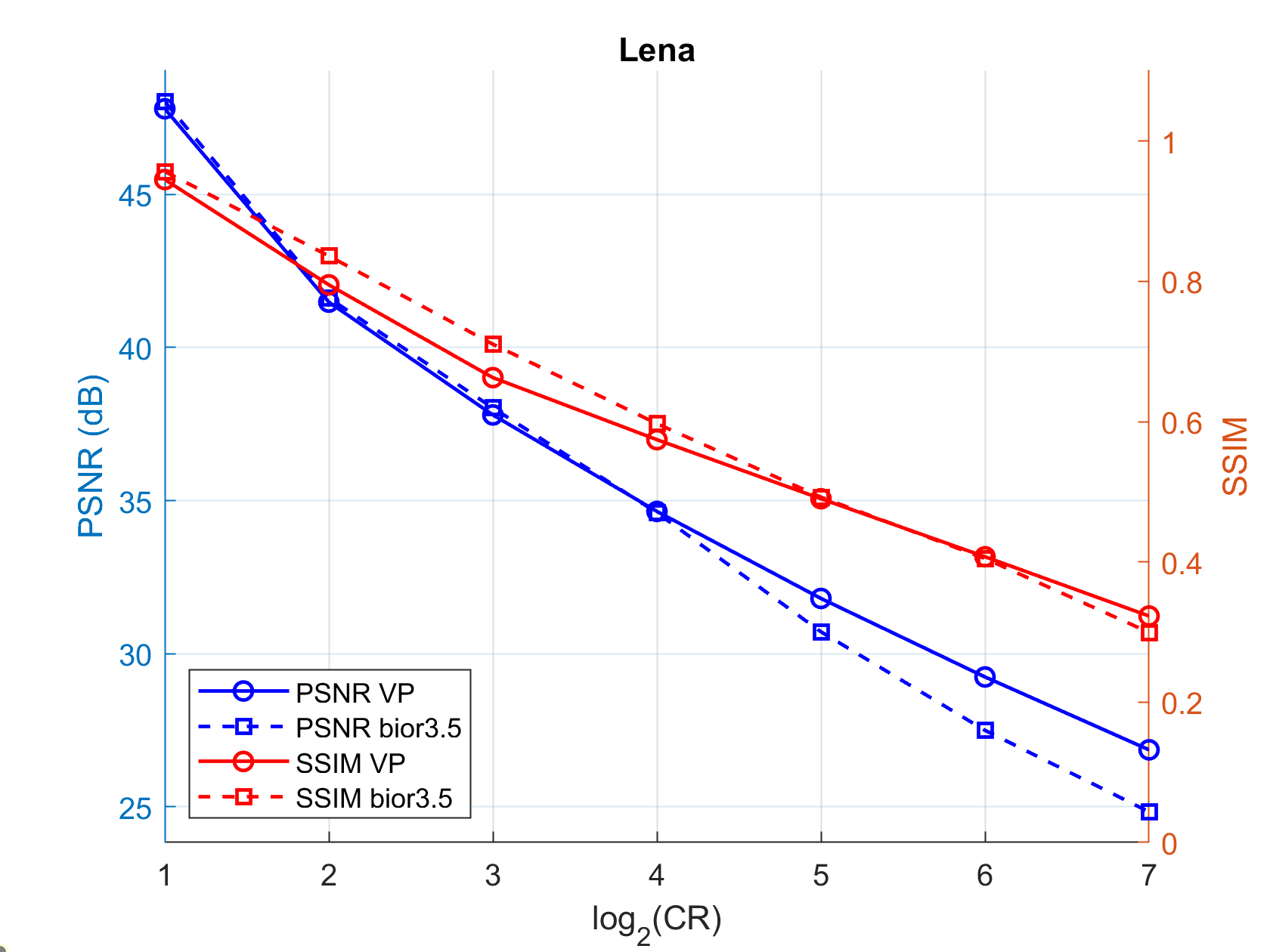}
 \includegraphics[width=6cm]{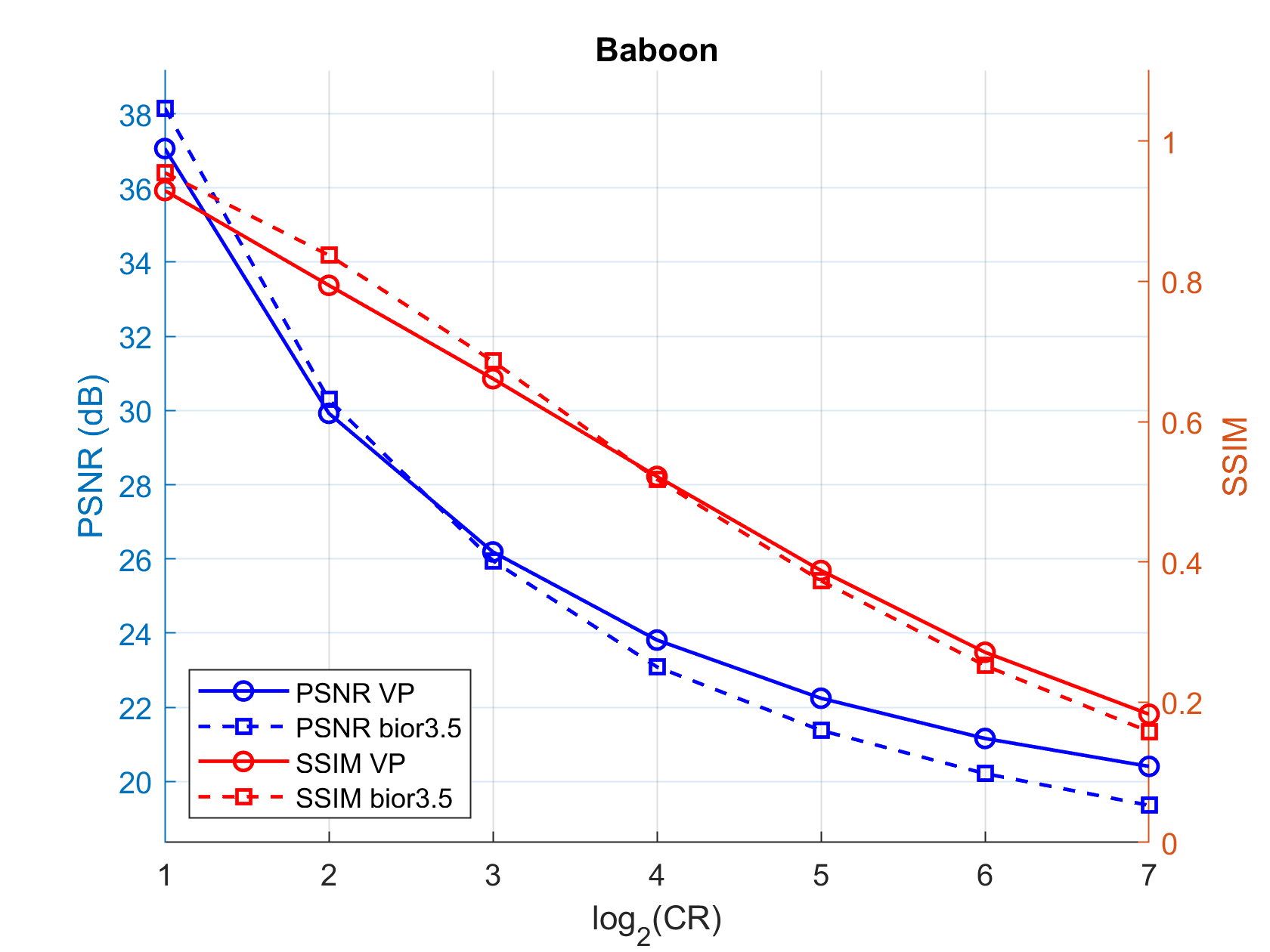}
 \includegraphics[width=6cm]{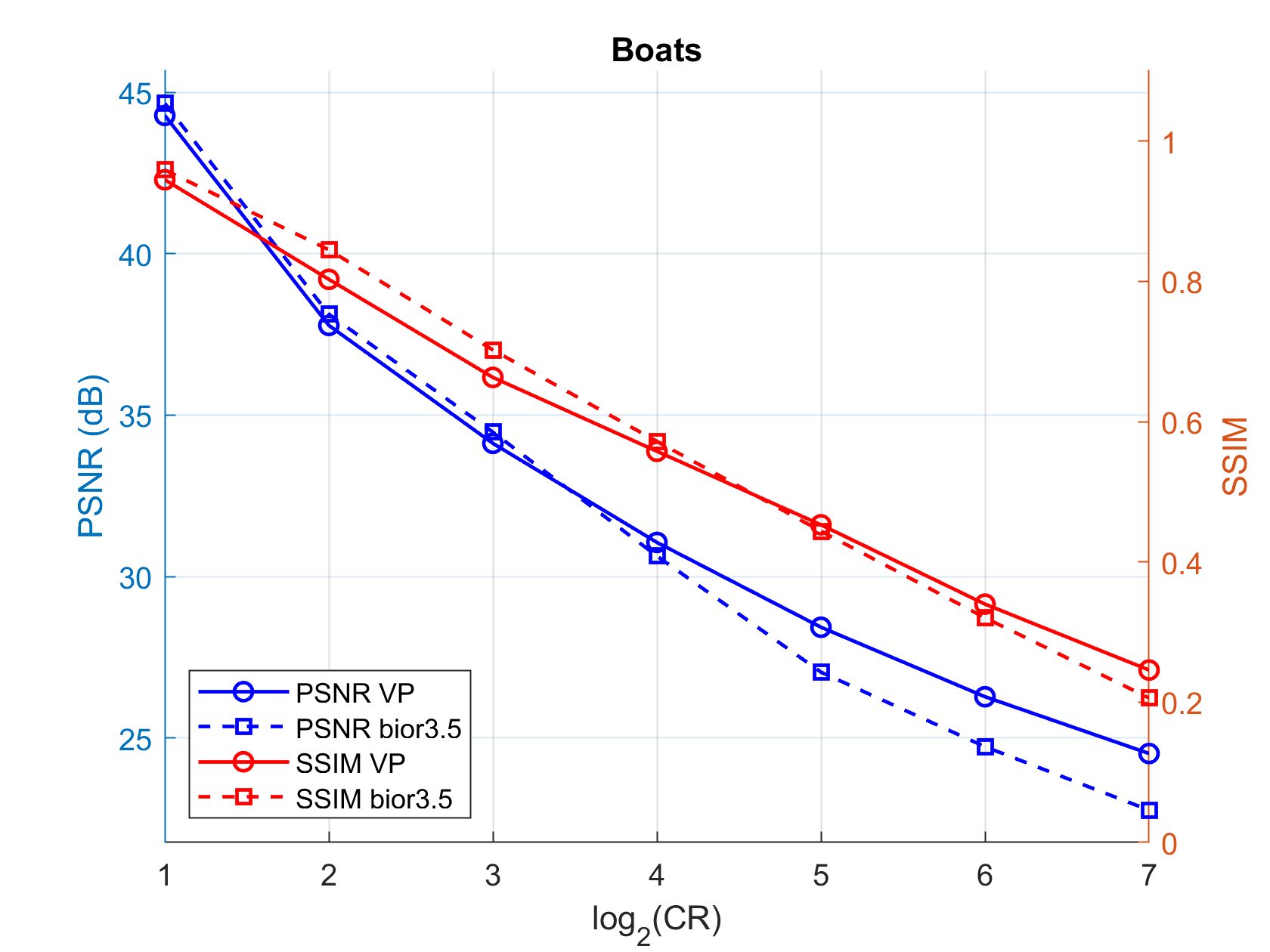}
 \includegraphics[width=6cm]{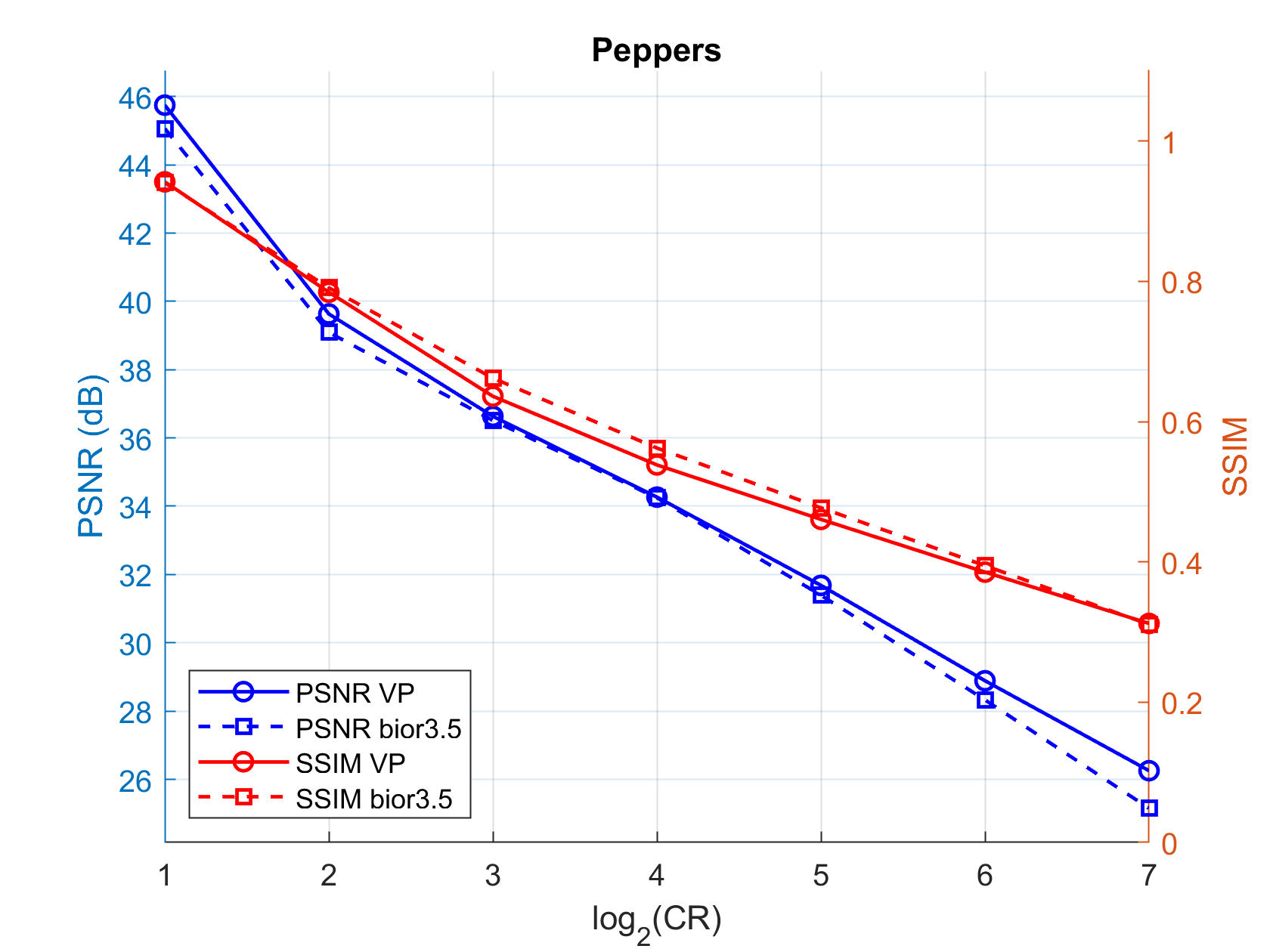}
\end{center}
\caption{Comparison of reconstruction quality between VP and Bior3.5 wavelet filters in terms of PSNR and SSIM as a function of $\log_2(\mbox{CR})$ over different test images \label{fig:Comparisons}}
\end{figure}
	
We now present a comparative analysis of the reconstruction quality between the VP basis and the biorthogonal spline wavelet filters bior3.5. We selected such filters  because they are widely used in image compression, thanks to their symmetry, compact support, and high number of vanishing moments. The evaluation is performed in terms of PSNR and SSIM, considering their variation as a function of $\log_2(\text{CR})$. For this comparison, we have fixed seven levels of compression ratio, namely $2^k$ for $k = 1, \dots, 7$. For a fair comparison, since the images have all a size of $512 \times 512$, we have considered 4 levels of decomposition for VP and 7 levels for the bior3.5 filters, so that, at the very last level, a similar amount of approximation coefficients are retained (4 and 6 coefficients, respectively). For the VP basis, a near-optimal $\theta$ parameter, providing nearly the best PSNR, was chosen from the interval $[0.5, 0.8]$ based on the experimental tests illustrated above.

From the plots in Fig. \ref{fig:Comparisons}, it can be observed that the VP basis generally yields comparable or slightly higher PSNR and SSIM values than the bior3.5 wavelets across the given range of compression ratios for all four test images. In particular, for highly textured images (e.g., Baboon), the VP basis often maintains better performance at moderate to higher compression ratios, suggesting it can handle  details effectively. For smoother images like Lena, its performance is still competitive and sometimes surpasses bior3.5. 

	\begin{figure}[t]
		\centering
		\includegraphics[width=5cm]{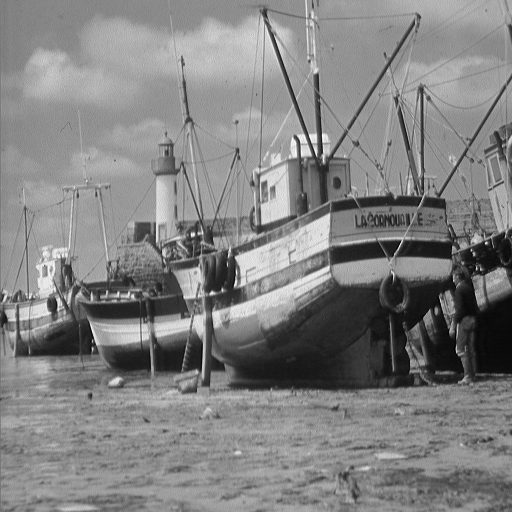} \\[0.3cm] 
		\includegraphics[width=10cm]{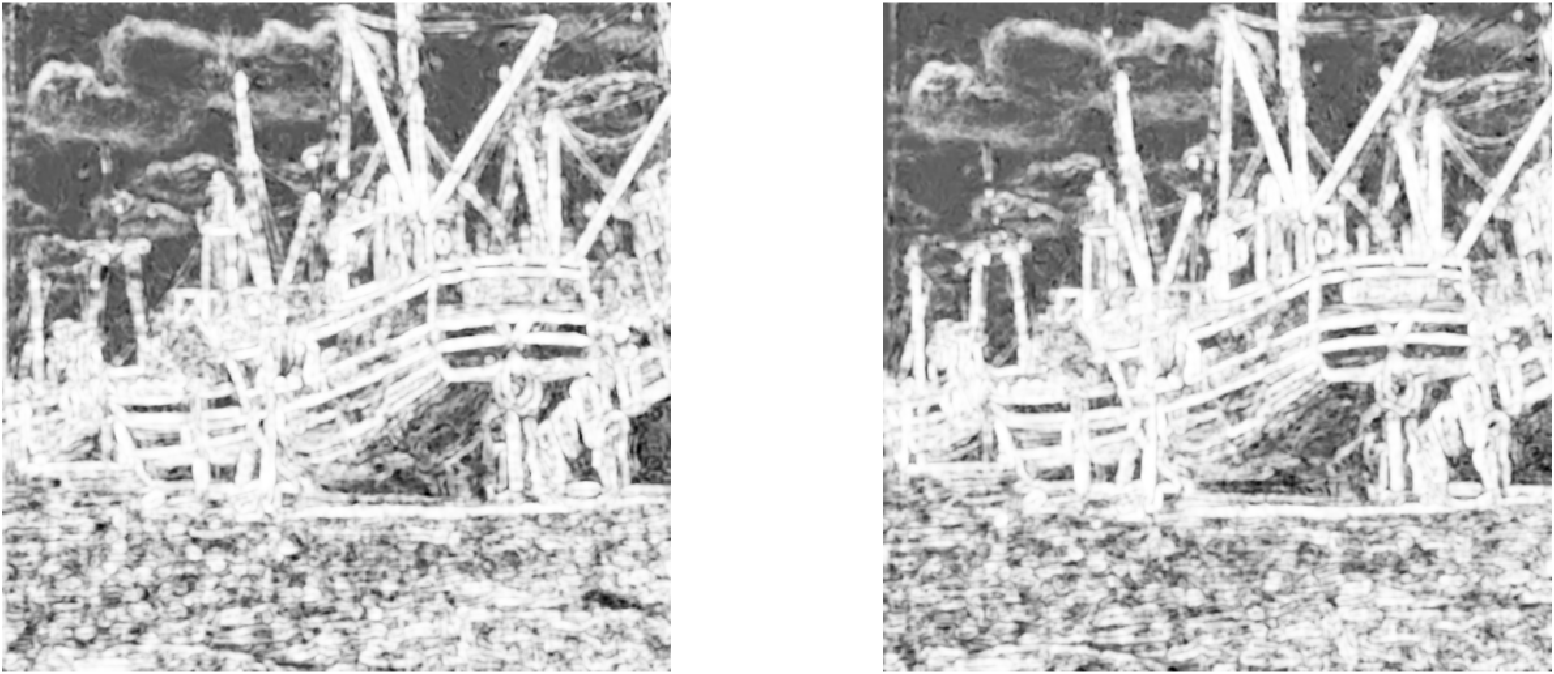} \\[0.3cm] 
		\includegraphics[width=10cm]{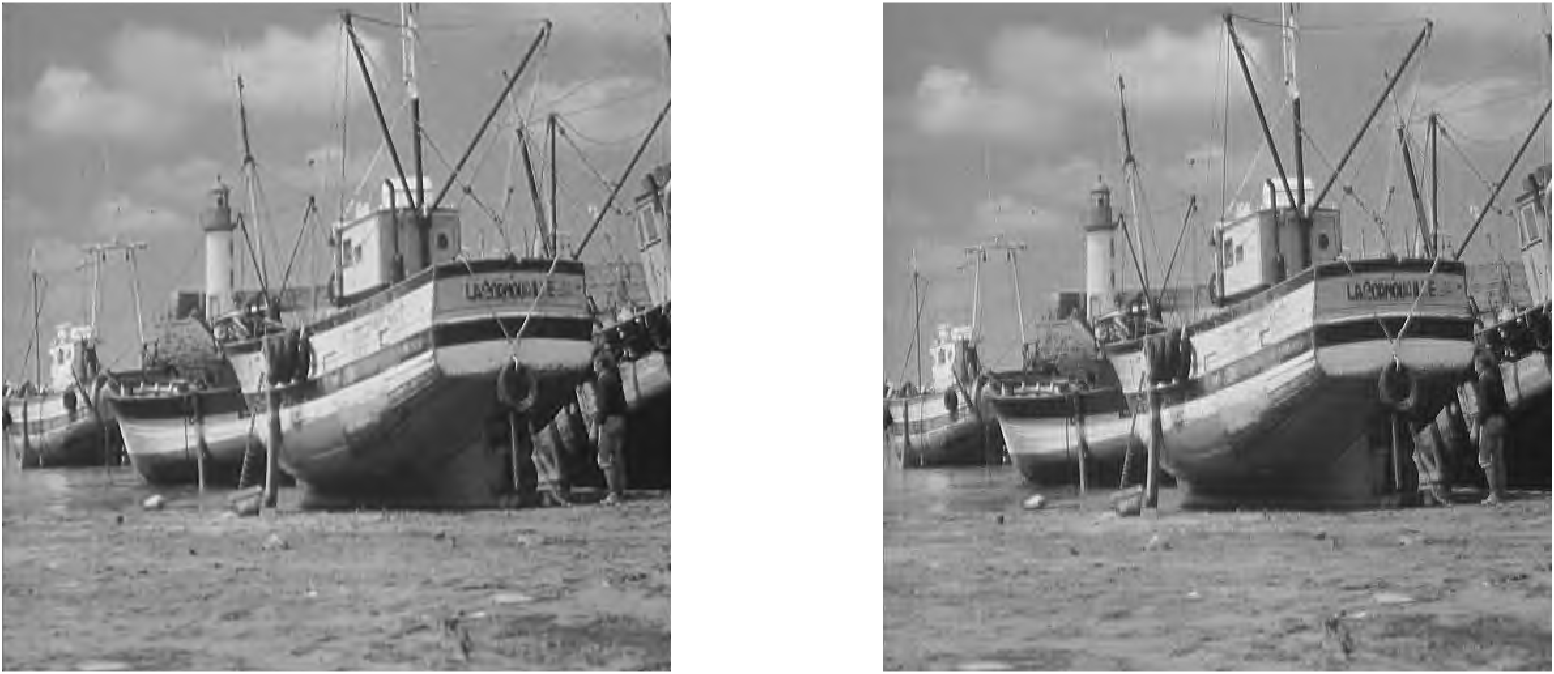} 
		\caption{Original boat image (top), full SSIM maps (middle), and reconstructed images (bottom) for the bior3.5 wavelet filter (left column) and the VP wavelet (right column).}
		\label{fig:boats_full}
	\end{figure}

	\begin{figure}[t]
		\centering
		\includegraphics[width=3cm, trim=120pt 200pt 120pt 0, clip]{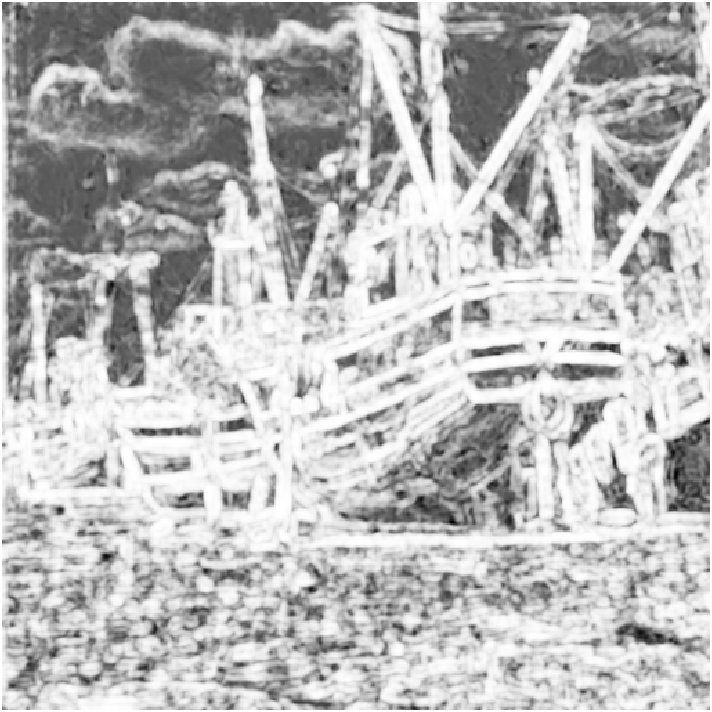} 
		\hspace{1cm}
		\includegraphics[width=3cm, trim=120pt 200pt 120pt 0, clip]{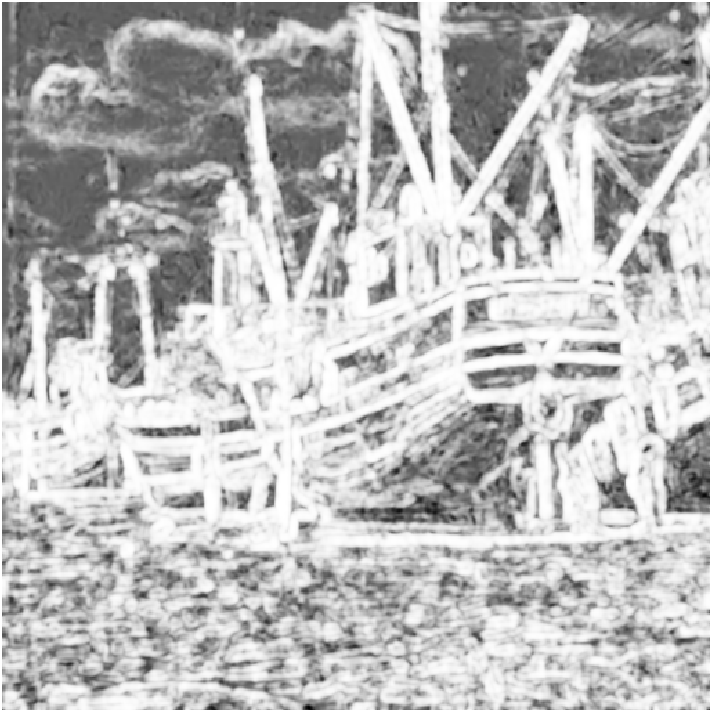} \\[0.3cm]
		\includegraphics[width=3cm, trim=120pt 200pt 120pt 0, clip]{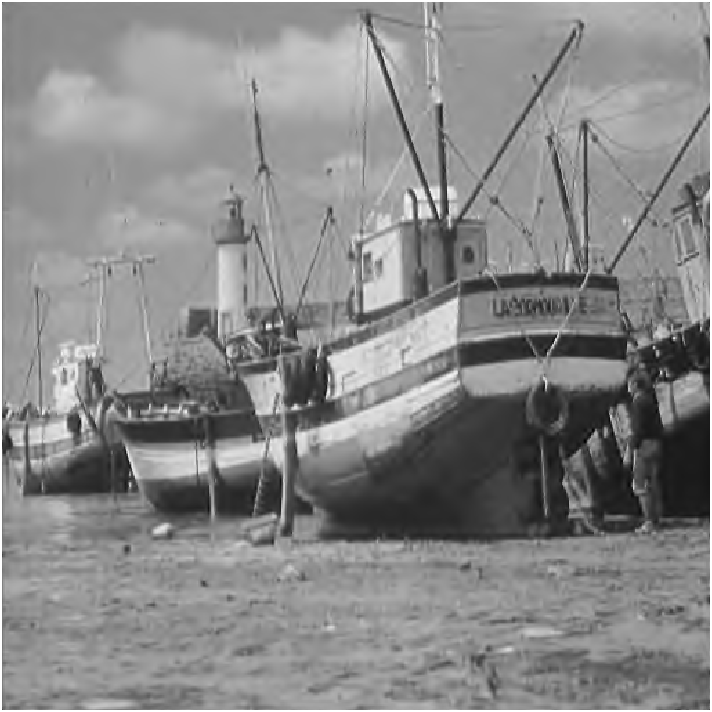} 
		\hspace{1cm}
		\includegraphics[width=3cm, trim=120pt 200pt 120pt 0, clip]{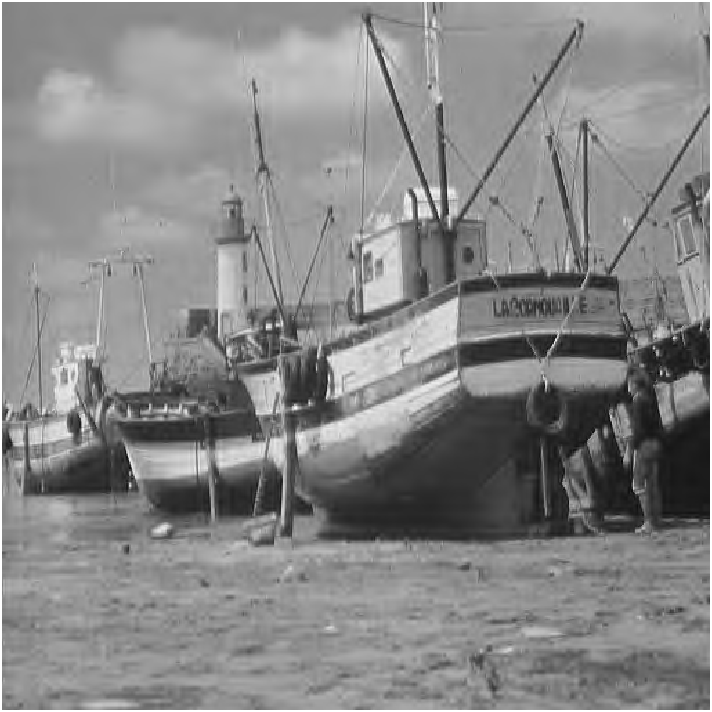}
		\caption{Zoomed views of the SSIM maps (top row) and the reconstructed images (bottom row) for the bior3.5 wavelet filter (left column) and the VP wavelet (right column).}
		\label{fig:boats_zoom}
	\end{figure}

\begin{figure}[t]
		\centering
				\includegraphics[width=7cm]{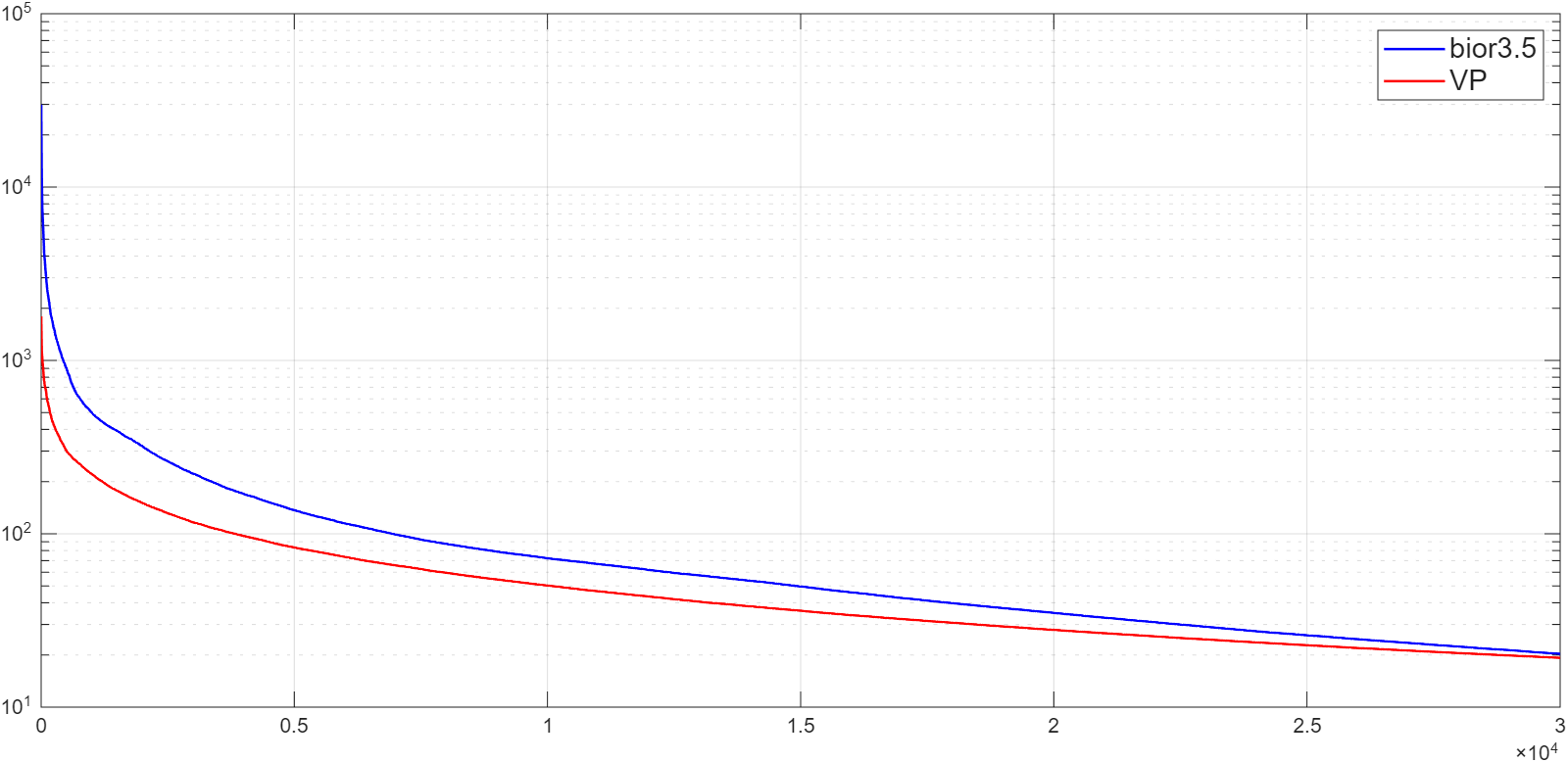}
		\caption{Absolute values (in log scale) of the first 30000 sorted wavelet detail coefficients for the biorthogonal (bor3.5) and VP decomposition applied to the Boats image}
		\label{fig:sorted_coeffs}
	\end{figure}

\MC{To illustrate the reconstruction quality visually, we show in Figs. \ref{fig:boats_full} and \ref{fig:boats_zoom} the results obtained for the "Boats" image after compression at $CR=2^{4}$,  retaining only about $6.25\%$ of the coefficients. The accompanying SSIM maps indicate the local structural similarity with the original image. Overall, the VP wavelets maintain sharper details. The zoomed views (Fig.~\ref{fig:boats_zoom}) focus on one of the ropes within the boat, highlighting that VP reconstructions exhibit better preservation of fine structures.
This visual evidence is consistent with the plot in Fig. \ref{fig:sorted_coeffs}, which shows, in logarithmic scale, the first 30000 wavelet coefficients (absolute values) sorted in decreasing order. The VP wavelet coefficients are smaller in magnitude than those of the classical transform, indicating that the VP decomposition captures most of the energy in the approximation subband and leaves only small residuals in the detail subband. }

\begin{figure}[t] 
\begin{center} 
\includegraphics[width=4cm]{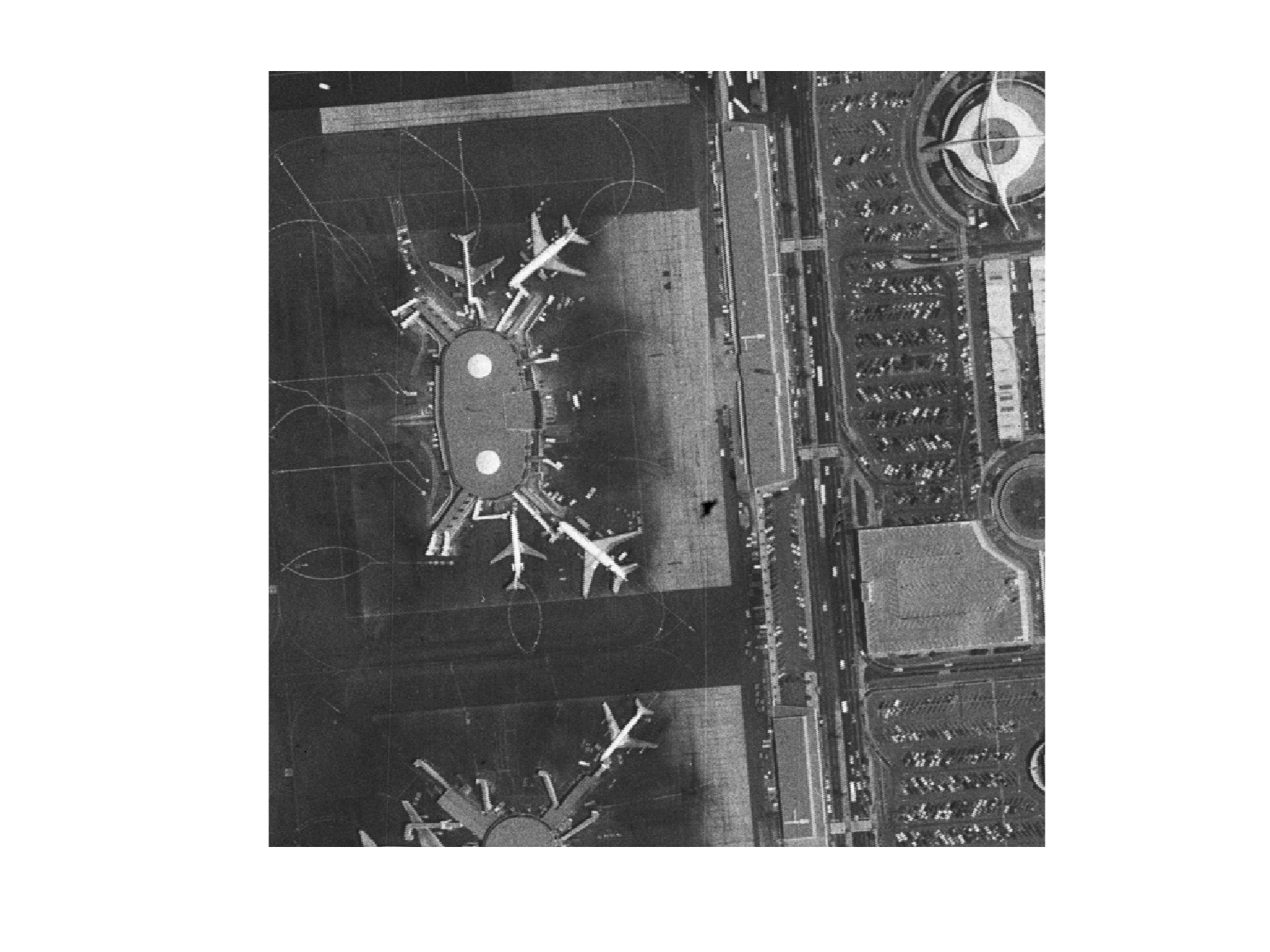}
\includegraphics[width=5cm]{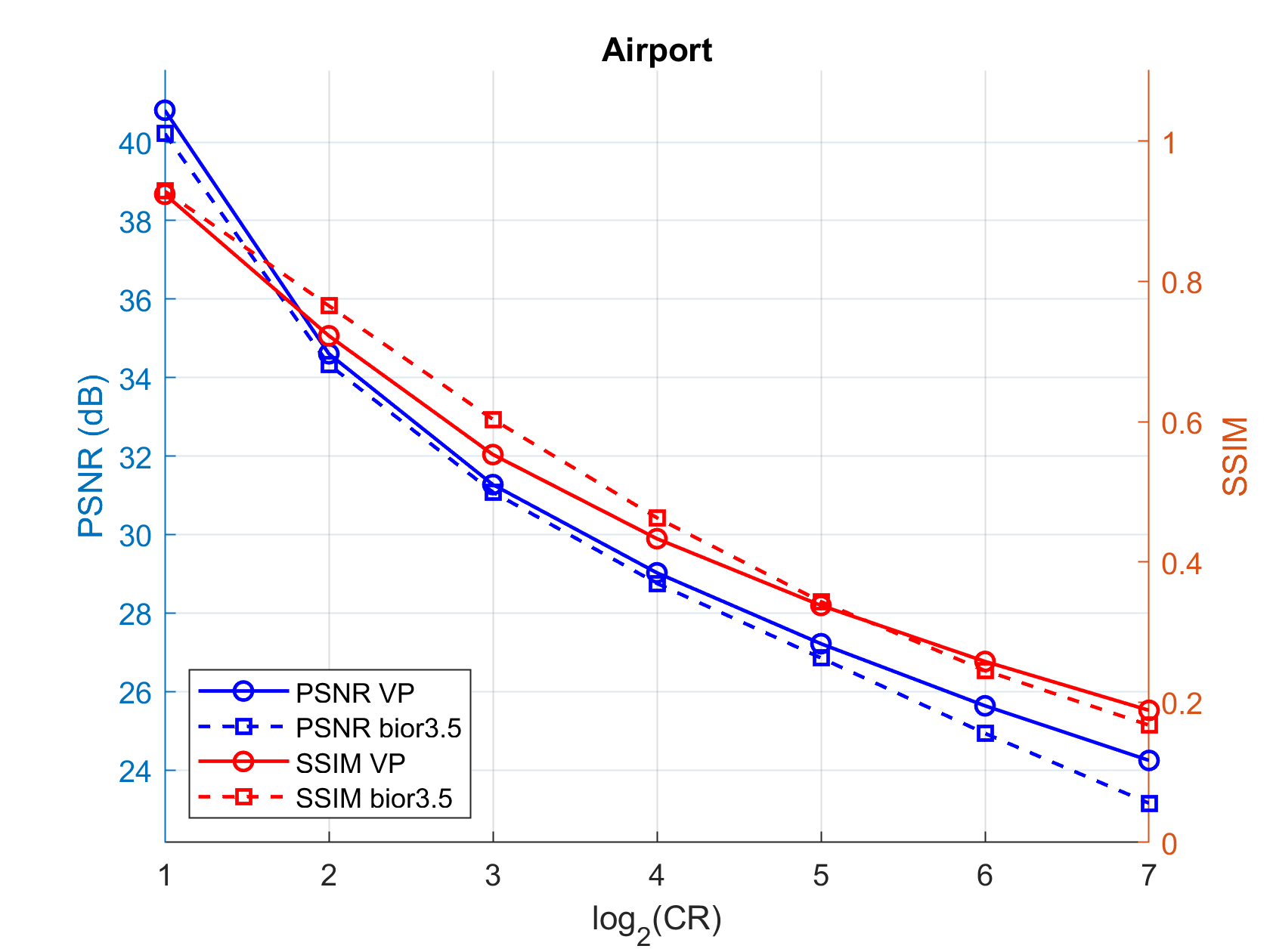} \\
\includegraphics[width=4cm]{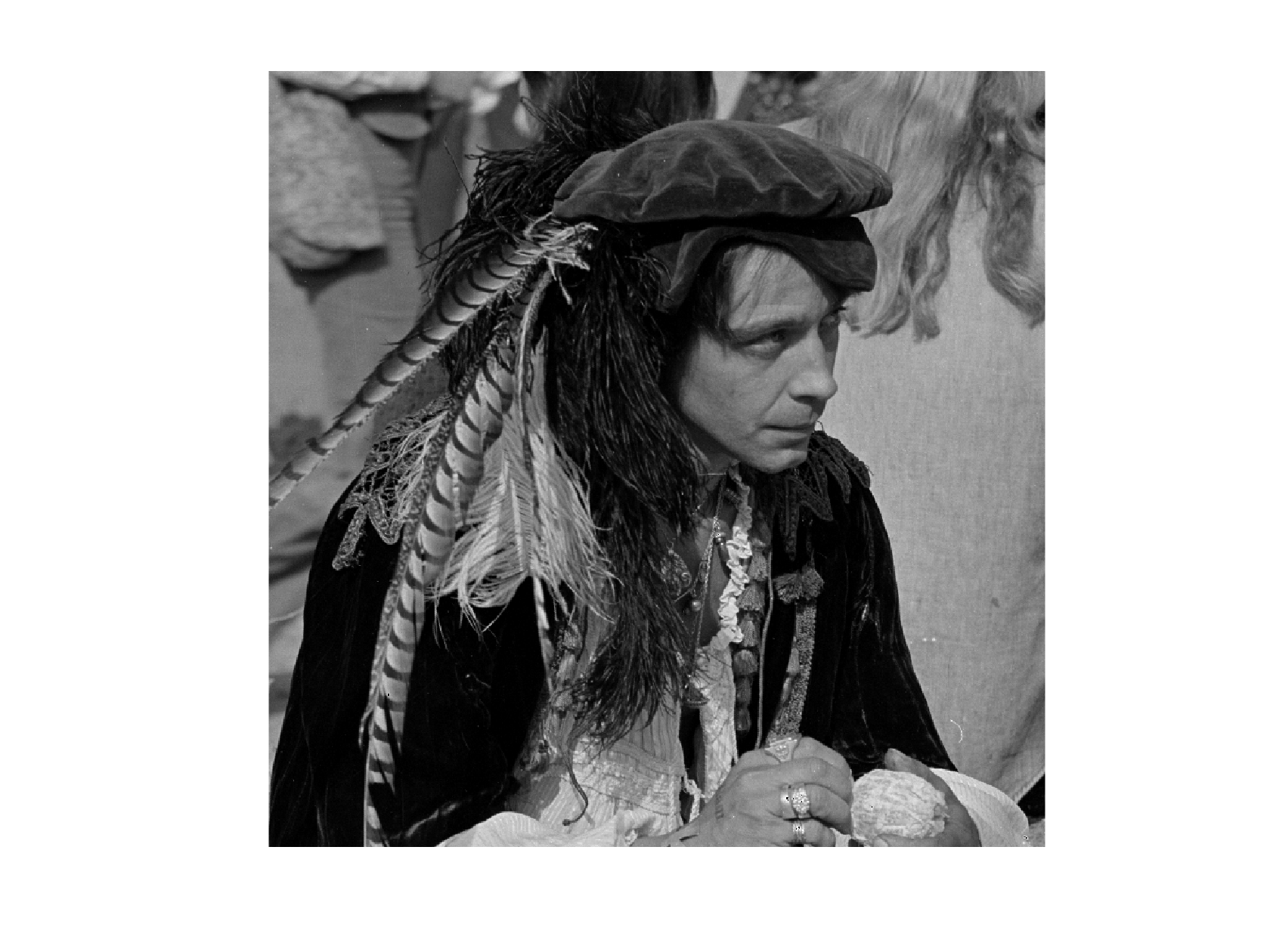}  
\includegraphics[width=5cm]{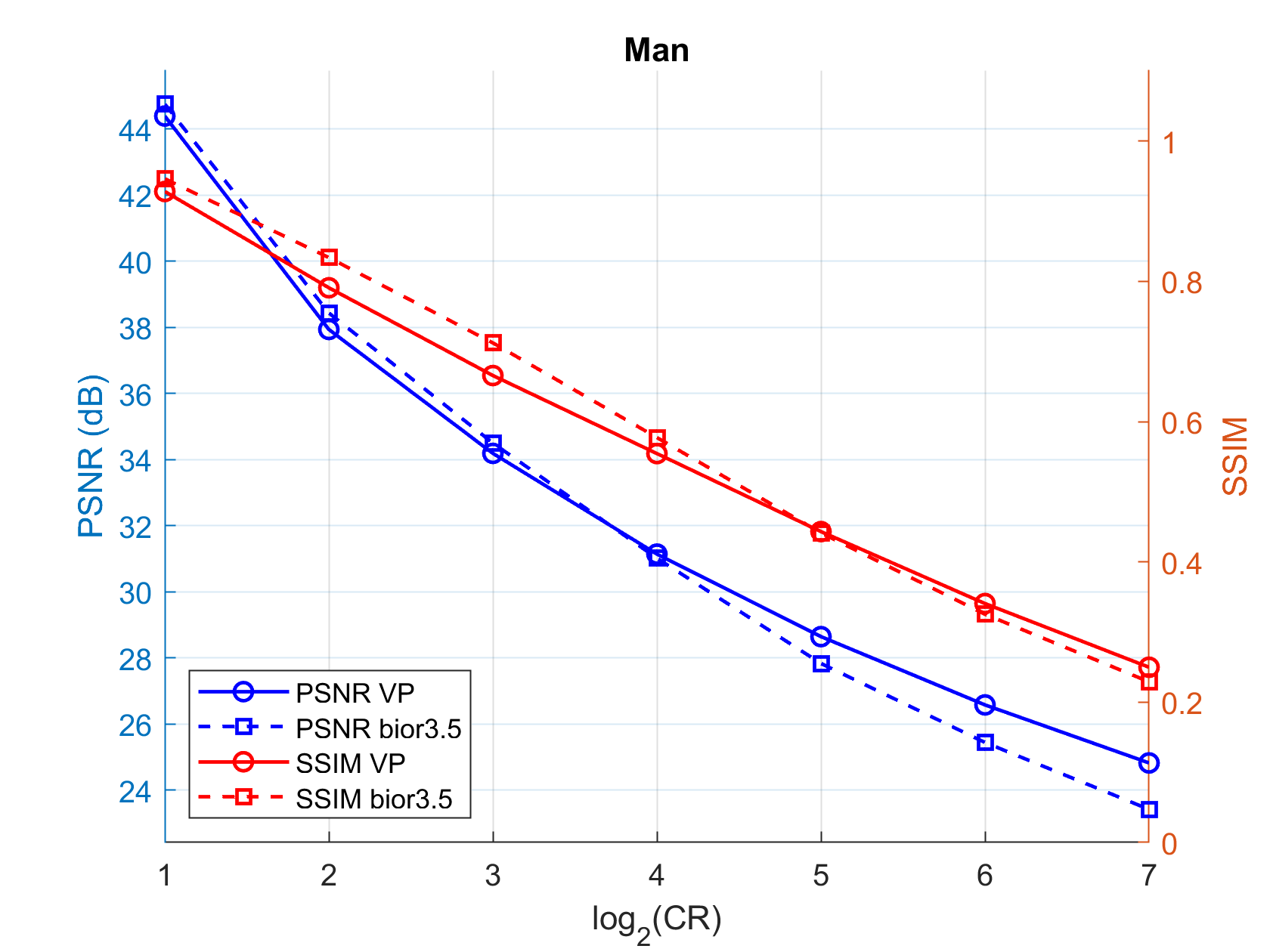} \caption{Comparison of reconstruction quality between VP and bior3.5 wavelet filters in terms
of PSNR and SSIM as a function of $\log_2(\mbox{CR})$ over  the two test images on the left image}\label{fig:test1024} 
\end{center}
\end{figure}

We now present results for two larger images, "Airport" and "Man" shown in Fig. \ref{fig:test1024} on the left. Originally these images have size $1024\times 1024$, but we considered a cropped version of them, fixing the size to $729\times 729$. This (power of 3)  size   ensures compatibility with our algorithm without requiring additional adjustments to enforce a power-of-three size at each decomposition step, as described in Section 3. 

\begin{figure}[t] 
\begin{center} 
\includegraphics[width=7cm]{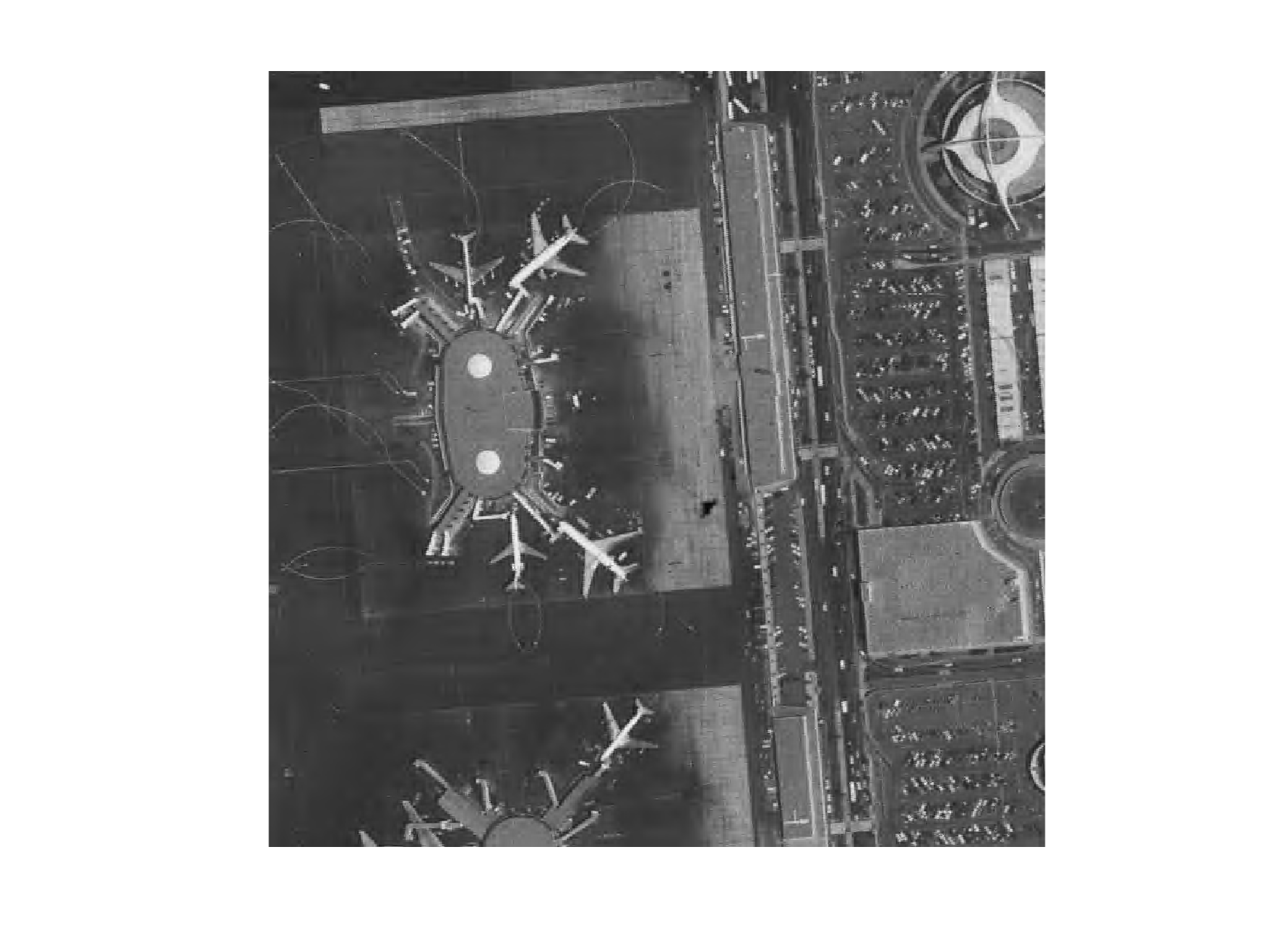} 
\includegraphics[width=7cm]{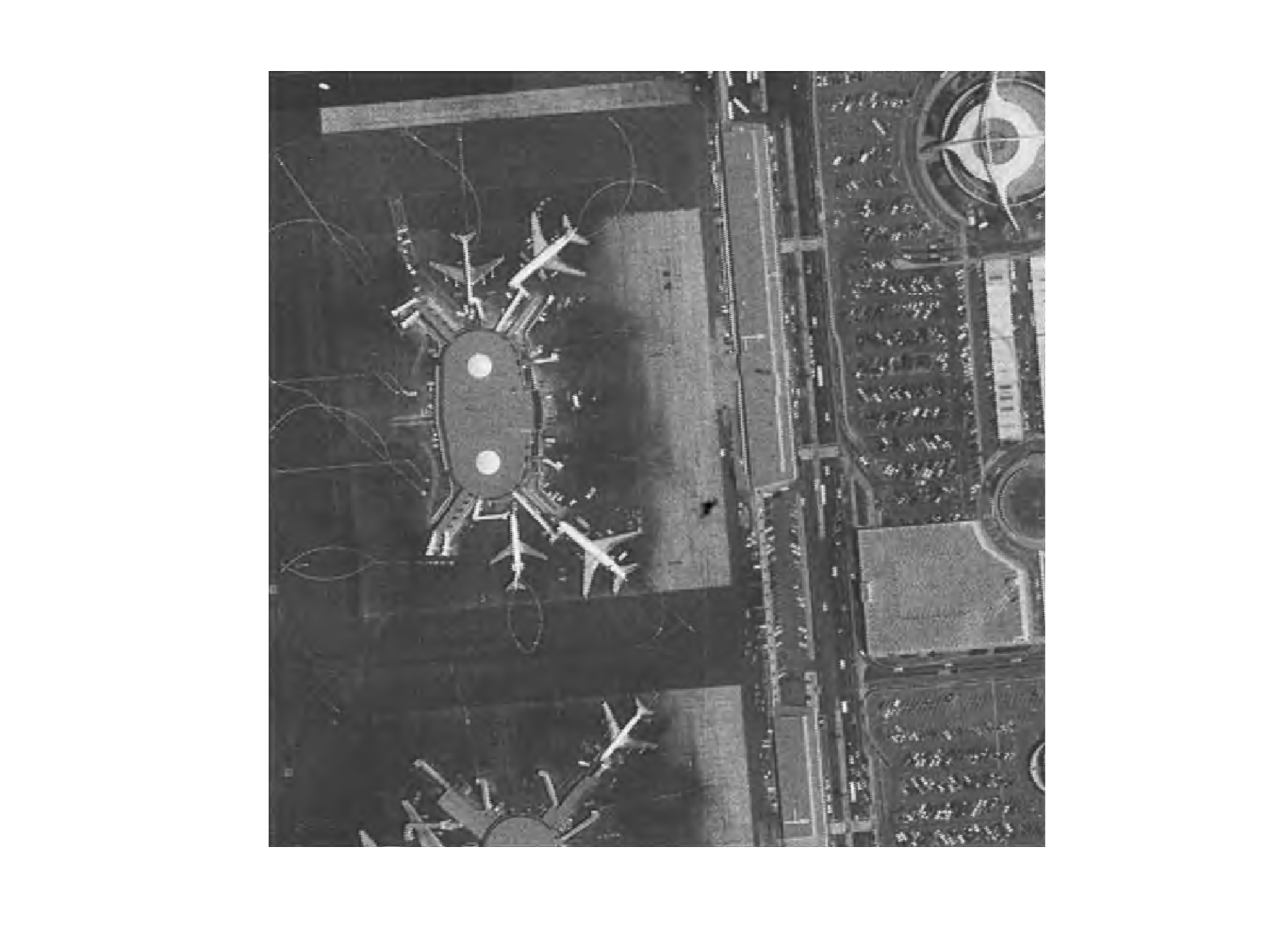} \caption{Reconstruction of the "Airport" image compressed  at CR $1/16$: VP (left) and bior3.5 (right).}\label{fig:rec} 
\end{center} 
\end{figure}

The results shown on the two graphs on the right of Fig. \ref{fig:test1024}  confirm that our basis outperforms bior3.5 for high-compression rates. While the performance appears similar or worse at lower compression rates, a closer visual inspection of the reconstructed images reveals that the VP basis produces quality, particularly in preserving fine details and reducing artifacts, as can be observed in Fig.\ref{fig:rec}  for the Airport image. 

As a further test, we compare the performance of VP wavelets with both the standard db2 and bior3.5 filters. For brevity, we show the results on three of the four 512 $\times$ 512 test images considered before, though  several other images with different sizes have been considered in our experiments displaying similar behavior.

In order to highlight that our decomposition enables a reduction in computational cost while achieving similar or higher PSNR and SSIM values, Table \ref{tab:1step} presents the results obtained after a single decomposition step for different percentages P of retained coefficients, i.e., keeping 50\%, 40\%, and 30\% of the image data. We note that stronger compression cannot be considered with just one decomposition step, as there would not be enough wavelet coefficients to compress. In this regard, we point out that, for instance, P= 25\%  is allowed for VP wavelets but not for db2 and bior3.5.

Finally, without any limit on the decomposition steps, we considered different percentages   P = 100/CR   of retained coefficients, including very large CR values. Table \ref{tab:multistep} displays the best results we obtained and the fewest decomposition steps required. As we can see, the performance of VP wavelets improves as   P   decreases (or   CR  increases), while for low CR values, db2 outperforms VP and bior3.5 wavelets, whose performances are almost comparable. In all cases, the decomposition steps required by VP wavelets are much fewer than those required by the other wavelets.

\begin{table*}[t]
\caption{Compression performance after just 1 decomposition step }
\label{tab:1step}
{\scriptsize
\begin{center}
\begin{tabular}{lc|rrc|rrc|rrc}
  &&\multicolumn{3}{c}{\small Lena} & \multicolumn{3}{|c}{\small Peppers} & \multicolumn{3}{|c}{\small Baboon}\\ \hline
&& PSNR & SSIM & Steps & PSNR  & SSIM & Steps & PSNR  & SSIM & Steps\\
\hline
\multirow{3}{*} {P=50\%}& VP &{\bf 47.354} &	{\bf 0.940}	&1 &{\bf 45.293}		&{\bf 0.936}	 &1&{\bf 36.581}		&0.924	 &1 \\
&db2 & 47.176	&0.937	&1 &45.046	&0.933	& 1 &36.466 &0.929	&1\\
&bior3.5& 46.241	&0.938	&1 &42.995 &0.908	&1 &36.128 &{\bf 0.939}	&1\\
\hline
\multirow{3}{*}{P=40\%}& VP  & {\bf 44.622} &	{\bf 0.893} & 1&{\bf 42.615}& {\bf 0.885} & 1&  33.527 &	0.880 &1\\
&db2&43.859	& 0.878&1 &  41.888 & 0.870 & 1& {\bf 33.648}& 0.887 & 1\\
&bior3.5 &43.859 & 0.878 &1 &40.464 & 0.846 & 1 & 32.687 &{\bf 0.899} & 1\\
\hline
\multirow{3}{*}{P=30\%}&VP & {\bf 42.049} &	{\bf 0.821}	&1 &{\bf 40.151}	 &{\bf 0.809} & 1 &{\bf 30.510}		&0.827	 &1 \\
&db2 & 39.726 &0.795	& 1& 38.388	&0.757	& 1 &28.109 &0.817	&1\\
&bior3.5& 40.362	&0.816	&1 &38.102	&0.759	&1 &28.526 &{\bf 0.833}	&1\\
\hline
\end{tabular}
\end{center}
}

\end{table*}
\begin{table*}[t]
\caption{Compression performance considering the best result achieved with the fewest decomposition steps  }
\label{tab:multistep}
{\scriptsize
\begin{center}
\begin{tabular}{lc|rrc|rrc|rrc}
  &&\multicolumn{3}{c}{\small Lena} & \multicolumn{3}{|c}{\small Peppers} & \multicolumn{3}{|c}{\small Baboon}\\ \hline
&& PSNR & SSIM & Steps& PSNR  & SSIM & Steps & PSNR  & SSIM & Steps\\
\hline
\multirow{3}{*}{P=50\%}&VP &47.797&	0.945	&3 &45.743		&0.942	 &3 &37.056	&0.930	 &3 \\
&db2 & {\bf 49.099}	&{\bf 0.956}	& 5 &{\bf 47.214}	&{\bf 0.957}	&5 &{\bf 38.538} &0.940	&5\\
&bior3.5& 47.349		&0.950	&5 &44.343 &0.931	& 4 &37.276 &{\bf 0.947}	&5\\
\hline
\multirow{3}{*}{P=25\%}& VP  & 41.481	&0.795 & 3 &
39.627 & 0.785 & 3&  29.920 &	0.795 &3\\
&db2&{\bf 42.218}	& 0.805& 6&  {\bf 40.640} & {\bf 0.815} & 5& {\bf 30.877}& {\bf 0.801} & 6\\
&bior3.5 &41.415 & {\bf 0.830} &5 &38.928 & 0.784 & 5 & 30.009 &0.830 & 6\\
\hline
\multirow{3}{*}{P=10\%}& VP & 36.730 &	0.632	&3 &35.860	 &0.600& 3&25.323		&0.619	 &3 \\
&db2 & 36.952 &0.623	& 6&{\bf 36.363}	&0.606	& 5&{\bf 25.726} &0.613	&6\\
&bior3.5& {\bf 37.022}		&{\bf 0.675}	&5 &35.889	&{\bf 0.631}	& 5&24.826 &{\bf 0.628}	&6\\
\hline
\multirow{3}{*}{P=5\%}& VP & 33.686&	0.546	&3 & 33.457 &	0.512	 &3 &23.347	&	{\bf 0.481} &3\\
&db2 &33.516 &0.526	&6 &{\bf 33.726}&	0.511& 8 &{\bf 23.460}&	0.470 &7 \\
&bior3.5 &{\bf 33.838}&	{\bf 0.578}	&6& 33.650&	{\bf 0.539} &5 &22.453	&	0.465& 6\\
\hline
\multirow{3}{*}{P=2\%}& VP & {\bf 30.138} &	0.440&3& {\bf 29.901}&	0.412	&3&21.502	&	{\bf 0.309} &4\\
&db2 & 29.585&	0.404	& 6&29.747 &0.405& 8 &{\bf 21.611} &	0.301 & 7\\
&bior3.5 & 29.900	&	{\bf 0.455} & 7 &29.698 &{\bf 0.433}& 6 &20.700&	0.294 & 6 \\
\hline
\multirow{3}{*}{P=1\%}& VP & {\bf 27.686} &	0.352&3 & {\bf 27.170}&	0.338	&3 &20.654	&	{\bf 0.212} &4 \\
&db2& 27.035	&	0.319	&8 &26.924 &0.330 &8 &{\bf 20.702}&	0.201&7\\
&bior3.5 &27.002	&	{\bf 0.359}& 6& 26.667& {\bf 0.351} &6 &19.880&	0.198 & 7\\
\hline
\multirow{3}{*}{P=0.5\%}& VP & {\bf 25.443} &	{\bf 0.275}&4& {\bf 24.666}&	0.267	&3&{\bf 20.041}	&	{\bf 0.144} &4 \\
&db2& 24.864	&	0.242	& 8&24.527 &0.256& 8&20.035&	0.131& 9\\
&bior3.5 &24.590	&	0.273& 7& 24.259&{\bf 0.277}& 9 &19.378&	0.129 & 7\\
\hline
\end{tabular}
\end{center}
}
\end{table*}

\section{Conclusion}
A new family of polynomial wavelets, based on de la Vallée Poussin interpolation and characterized by a free parameter, was recently introduced in \cite{TheVanBar2024}, providing a solid theoretical foundation for future applications. In this study, after detailing the corresponding transform implementation, we evaluated the performance of these VP wavelets in two typical tasks — signal denoising and image compression — comparing them with classical wavelets.

The experimental results confirmed that VP wavelets offer competitive and satisfactory performance in both tasks, showing that VP wavelets are a viable alternative to standard wavelet approaches, providing robust reconstruction quality across different types of data. 

In the experiments, the role of the free parameter  $ \theta$  was also observed experimentally.

For signal denoising, the value  $ \theta = 0.1 $ was found to give better performance for the considered test signals. With this choice, the VP wavelet performs significantly better in the case of some test signals ("quadchirp" in particular), and for small input SNR, the VP wavelet gives results comparable to those of the standard wavelet approaches.
    
In image compression, we observed that the value of  $ \theta$  providing the best performance is strongly influenced by the image data and the fixed compression ratio. For the ``best'' parameter choice, all the experiments revealed much better performance for larger compression ratios (i.e., smaller percentages of retained data) in comparison with Daubechies and biorthogonal wavelets.

Nevertheless, with only one decomposition step, VP wavelets always outperform both Daubechies and biorthogonal wavelets.
 Indeed, compared to classical wavelets, VP wavelets have the advantage of considering multiples of 3 instead of 2 in the downsampling step of a multiresolution analysis. 
   
   As a further advantage, they do not require periodization or zero-padding techniques because they are defined on a compact interval.

A disadvantage is their non-orthogonality, which led us to introduce appropriate normalization factors for energy preservation. These factors have been determined theoretically, starting from the Riesz stability properties in \eqref{Riesz-sca} and \eqref{Riesz-wav}, which establish an equivalence between continuous and discrete norms. In such equivalences, the involved constants are independent of the data and the resolution level but depend on the $ \theta$ parameter. For denoising, these normalization factors have been further refined, with experimental values obtained.

In future research, we aim to find better values for the normalization factors. We have observed that these factors are dependent on the choice of  $ \theta$, which deserves further investigation.

We also aim to explore other applications that can best exploit the regularity, interpolation, and uniform near-best approximation properties offered by VP wavelets.

\section*{Acknowledgments}
The first two authors are members of RITA (Research Italian Network on Approximation),
of UMI-TAA research group, and Gruppo Nazionale Calcolo Scientifico-Istituto Nazionale
di Alta Matematica (GNCS-INdAM) which partially supported this work. The first author was also partially supported by the EU under the Italian National Recovery and Resilience Plan (PNRR) of NextGenerationEU, partnership on “Telecommunications of the Future” (PE00000001 - program “RESTART”).
The research of the third author was partially supported by the Fund for Scientific Research -- Flanders (Belgium), project G0B0123N.

\MC{The authors thank the anonymous referees for their helpful comments and suggestions, which improved the quality of the paper.}


\end{document}